\documentclass[a4paper]{article}
\usepackage[margin=0.8in,footskip=0.25in]{geometry}
\usepackage{amsmath, amsthm, amscd, amsfonts, amssymb, graphicx, color,float,pgf,tikz}
\usepackage[russian,english]{babel}
\usepackage{amssymb,fontenc}
\usepackage{cite}
\usepackage{latexsym,wasysym,mathrsfs}
\usepackage{hyperref}
\usetikzlibrary{arrows}

\theoremstyle{remark}

\numberwithin{equation}{section}

\begin{document}

\begin{center}

\textbf{\LARGE{Necessary conditions for a minimum in classical calculus of variation problems in the presence of various degenerations}}
\end{center}

 \begin{center}
 \textbf{M. J. Mardanov }${\boldsymbol{\mathrm{\ }}}^{\boldsymbol{\mathrm{1}}\boldsymbol{\mathrm{,\ }}}$\textbf{, T. K. Melikov }${\boldsymbol{\mathrm{\ }}}^{\boldsymbol{\mathrm{2}}}$\textbf{, S. T. Melik }${\boldsymbol{\mathrm{\ }}}^{\boldsymbol{\mathrm{3}}}$
 \end{center}

 \begin{center}
 ${\ }^1$\textit{IMM of ANAS, Vahabzadeh str. 9, Az 1141, Baku; Baku State University, Baku 1148, Z.Khalilov str. 23.}

 \end{center}
\begin{center}
 \noindent E-mail: misirmardanov@yahoo.com 

 $^2$\textit{ IMM of ANAS, Vahabzadeh str. 9, Az 1141, Baku; Institute of Control Systems of ANAS, Baku, Azerbaijan.}

\noindent E-mail: t.melik@rambler.ru 

\noindent  $^3$\textit{ Baku Higher Oil School, Baku, Azerbaijan; IMM of ANAS, Vahabzadeh str. 9, Az 1141, Baku;. }

\noindent E-mail: saminmelik@gmail.com
 \end{center} 

\textbf{Abstract. }In the paper, we offer a method for studying an extremal in the classical calculus of variation in the presence of various degenerations. This method is based on introduction of Weierstrass type variations characterized by a numerical parameter. To obtain more effective results, introduced variations are used in two forms: in the form of variations on the right with respect to the given point, and in the form of variations on the left with respect to the same point. The research is conducted under the assumption that along the considered extremal the Weierstrass condition and the Legendre condition degenerate, i.e. they are fulfilled as equalities at separate points or on some intervals. Two types of new necessary conditions are obtained: of equality type and of inequality type conditions for a strong and also a weak local minimum. Given specific examples and counterexample show that some of the necessary minimum conditions obtained in this article are strengthening and refining  of the corresponding known results in this direction.

\textbf{Keywords:} calculus of variation, strong (weak) local minimum, necessary condition type equality (inequality), extremal, degeneration at the point (on the interval)

\section{Introduction and Problem Statement} 
We consider the following vector problem of the classical calculus of variations:    
\begin{equation} \label{GrindEQ__1_1_} 
J\left(x\left(\cdot \right)\right)=\int _{t_{0} }^{t_{1} }L\left(t,\, x\left(t\right),\, \, \dot{x}\left(t\right)\right)\,  dt\to {\mathop{\min }\limits_{{}_{x(\cdot )} }} ,    
\end{equation} 
\begin{equation} \label{GrindEQ__1_2_} 
x\left(t_{0} \right)=x_{0} ,  x\left(t_{1} \right)=x_{1} , \,\,\,  x_{0} , x_{1} \in R^{n} ,    
\end{equation} 
where $R^{n} $ is $n$ -dimensional Euclidean space, $L\, \left(\cdot \right)$ is a given function, $x_{0} ,\, x_{1} ,\, t_{0} ,\, t_{1} $ are the given points. The function $L\left(t,x,y\right):\left(a,b\right)\times R^{n} \times R^{n} \to R:=\left(-\infty ,\, +\infty \right)$, called an integrant, is assumed to be continuously differentiable by a totality of variables, where $\left(a,\, b\right)$ is an interval, and $\left[t_{0} ,\, t_{1} \right]\subset \left(a,b\right)$. The sought-for function $x\left(\cdot \right):\left[t_{0} ,\, t_{1} \right]=:I\to R^{n} $ is a piecewise-smooth vector-function, i.e. is continuous, and its derivative is continuous everywhere on $I$, except for a finitely many points $\tau _{i} \in \left(t_{0} ,\, t_{1} \right),\, \, i=\overline{1,m}$, and at the points $\tau _{i} $ the derivative function $\dot{x}\left(\cdot \right)$ has first kind discontinuities (at the points $t_{0} $and $t_{1} $the values of the derivative function $\dot{x}\left(\cdot \right)$ are finite on the right and left, respectively). Herewith, as a rule, the points $\tau _{i} ,\, \, \, i=\overline{1,m}$ are called angular points for the functions $x\left(\cdot \right)$. We denote the set of all piecewise-smooth functions on $\left[t_{0} ,\, t_{1} \right]$ by $KC^{1} \left(I,\, R^{n} \right)$.

We call the functions $x\left(\cdot \right)\in KC^{1} \left(I,\, R^{n} \right)$ satisfying the boundary condition \eqref{GrindEQ__1_2_}, admissible functions. Obviously, if $\bar{x}\, \left(\cdot \right)$ is a fixed admissible function, then for every $\theta \in \left[t_{0} ,\, t_{1} \right)\, \, \left(\theta \in \left(t_{0} ,\, t_{1} \right]\, \right)$ there exists a number $\alpha >0$ such that $\bar{x}\, \left(\cdot \right)$ is continuously differentiable on the semi-interval $\left[\theta ,\, \theta +\alpha \right)\subset I\, \left(\left(\theta -\alpha ,\, \theta \right]\subset I\right)$. We call this statement the property $P\left(\theta +;\, \bar{x}\left(\cdot \right),\, \, \alpha \right)$ $\left(P\left(\theta -;\, \bar{x}\left(\cdot \right),\, \, \alpha \right)\right)$ for the function $\bar{x}\, \left(\cdot \right)$, and we will use it in the future.

Development of theory of classical variational calculus was set out, for example, in the papers and monographes [1, 4-6, 11, 15-17, 29, 30, 33 and etc.], where detailed review of main results obtained for problem (1.1), (1.2) and their various essential generalizations were given.

Recall (see [11, p. ]) some notions from the classical calculus of variations. The admissible function $\bar{x}\, \left(\cdot \right)$ is said to be a strong (weak) local minimum in problem \eqref{GrindEQ__1_1_}, \eqref{GrindEQ__1_2_} if there exists such a number $\bar{\delta }>0\, \, \, \left(\hat{\delta }>0\right)$  that the inequality$J\left(x\left(\cdot \right)\right)\ge J\left(\bar{x}\left(\cdot \right)\right)$ is fulfilled for all admissible functions $x\left(\cdot \right)$ for which 
\[\left\| \, x\left(\cdot \right)-\bar{x}\left(\cdot \right)\, \right\| _{C\left(I,R^{n} \right)} \le \bar{\delta }\, \, \, \, \left(\max \, \left\{\, \left\| x\, \left(\cdot \right)-\hat{x}\, \left(\cdot \right)\right\| _{C\left(I,R^{n} \right)} ,\, \left\| \dot{x}\, \left(\cdot \right)-\dot{\bar{x}}\, \left(\cdot \right)\, \right\| _{L_{\infty } \left(I,R^{n} \right)} \right\}\le \hat{\delta }\right).  \] 

In this case we say that the admissible function $\bar{x}\, \left(\cdot \right)$ affords a strong (weak) local minimum in problem \eqref{GrindEQ__1_1_}, \eqref{GrindEQ__1_2_} with $\bar{\delta }\, \, \left(\hat{\delta }\right)$ -neighborhood. Obviously, any strong local minimum at the same time is weak as well, but the opposite is not always true (see \cite{20}).

We also recall (see e.g. \cite{4}) some known necessary conditions for a strong and weak local minimum for the considered problem \eqref{GrindEQ__1_1_}, \eqref{GrindEQ__1_2_}. Let $\{\tau\}\subset (t_{0}, t_1) $  be a set of angular points of the admissible function $\overline{x}(\cdot)$. Then:  

(i) if $\bar{x}\left(\cdot \right)$ is a weak local minimum in problem \eqref{GrindEQ__1_1_}, \eqref{GrindEQ__1_2_}, then at the points $t\in I\setminus \{\tau\}$ it satisfies the Euler equation, i.e. for every $t\in I\setminus \{\tau\}$ we have the equality
\begin{equation} \label{GrindEQ__1_3_} 
\frac{d}{dt} L_{\dot{x}} \left(t,\, \bar{x}\, \left(t\right),\, \dot{\bar{x}}\left(t\right)\right)=L_{x} \left(t,\, \bar{x}\, \left(t\right),\, \dot{\bar{x}}\left(t\right)\right),    
\end{equation} 
and also for every $t\in \{\tau\}$ along the function $\overline{x}(\cdot)$ the Weierstrass-Erdman condition is fulfilled, i.e. the following equalities are valid
\begin{equation} \label{GrindEQ__1_4_} 
\bar{L}_{x} \left(t-\right)=\bar{L}_{x} \left(t+\right),\, \, \bar{L}\left(t-\right)-\dot{\bar{x}}\left(t-\right)\, \bar{L}_{\dot{x}} \left(t-\right)=\bar{L}\left(t+\right)-\dot{\bar{x}}\left(t+\right)\bar{L}_{\dot{x}} \left(t+0\right);  
\end{equation} 

(ii) if  $\bar{x}\left(\cdot \right)$ is a strong local minimum in problem \eqref{GrindEQ__1_1_}, \eqref{GrindEQ__1_2_}, then the Weierstrass condition is fulfilled along it, i.e. for all  $\xi \in R^n$ the following inequalities are valid:
\[E\left(\bar{L}\right)\, \left(t,\, \xi \right)\ge 0,\, \, \, \forall t\in \left(t_{0} ,\, t_{1} \right)\backslash \left\{\tau \right\},\] 
\begin{equation} \label{GrindEQ__1_5_} 
E \left(\bar{L}\right)\, \left(t\pm ,\, \xi \right)\ge 0,\, \, \, \forall t\in \left\{\tau \right\},\, \, \, E\left(\bar{L}\right)\left(t_{0} +,\xi \right)\ge 0,\, \, \, E\left(\bar{L}\right)\, \left(t_{1} -,\, \xi \right)\ge 0. 
\end{equation} 

Here, $\bar{L}\left(t\right):=L\left(t,\, \bar{x}\left(t\right),\, \, \dot{\bar{x}}\left(t\right)\right),$$\bar{L}_{y} \left(t\right):=\bar{L}_{y} \left(t,\, \bar{x}\left(t\right),\, \, \dot{\bar{x}}\left(t\right)\right),\, \, y\in \left\{x,\, \dot{x}\right\}$,   
\begin{equation} \label{GrindEQ__1_6_} 
E\left(\bar{L}\right)\, \left(t,\, \xi \right):=E\left(L\right)\, \left(t,\, \bar{x}\left(t\right),\, \dot{\bar{x}}\, \left(t\right),\, \dot{\bar{x}}\, \left(t\right)+\xi \right)=L\left(t,\, \bar{x}\left(t\right),\, \dot{\bar{x}}\, \left(t\right)+\xi \, \right)-\, \bar{L}\left(t\right)\, -\bar{L}_{\dot{x}}^{\, T} \left(t\right)\, \xi , 
\end{equation} 
where the symbol $T$ denotes the transposition operation and $E\left(L\right)\left(t,\, x,\, y,\, z\right)$ is a Weierstrass  function for problem \eqref{GrindEQ__1_1_}, \eqref{GrindEQ__1_2_} and is determined in the form: 
\begin{equation} \label{GrindEQ__1_7_} 
E\left(L\right)\, \left(t,x,y,z\right)=L\left(t,x,z\right)-L\left(t,x,y\right)-L_{y}^{T} \left(t,x,y\right)\left(z-y\right).   
\end{equation} 

Underline that here and in what follows the symbol $F\left(t+\right)$ $\left(F\left(t-\right)\right)$ means a right (left) hand limit of the function $F\left(\cdot \right)$ at the point $t$, furthermore, fulfillment of equality \eqref{GrindEQ__1_3_} for $t=t_{0} $ $\left(t=t_{1} \right)$ is understood as a right (left) hand limit at the point $t_{0} $ $\left(t_{1} \right)$.

Following \cite{3, 21, 22}, we give a local modification of necessary condition for a minimum \eqref{GrindEQ__1_5_}. Let $\bar{x}\, \left(\cdot \right)$ be a weak local minimum in problem \eqref{GrindEQ__1_1_}, \eqref{GrindEQ__1_2_}. Then there exists a number $\delta >0$ at which the following inequalities are fulfilled:
\begin{equation} \label{GrindEQ__1_8_} 
\begin{array}{l} {E\left(\bar{L}\right)\, \left(t,\, \xi \right)\ge 0,\, \, \, \, \forall t\in \left(t_{0} ,\, t_{1} \right)\backslash \left\{\tau \right\},\, \, \, \, E\left(\bar{L}\right)\, \left(t\pm ,\xi \right)\ge 0,\, \, \forall t\in \left\{\tau \right\},\, } \\ {E \left(\bar{L}\right)\, \left(t_{0} +,\, \xi \right)\ge 0,\, \, \, E \left(\bar{L}\right)\, \left(t_{1} -,\, \xi \right)\ge 0,\, \, \, \, \forall \xi \in B_{\delta } \left(0\right).} \end{array} 
\end{equation} 
Here and in what follows, the symbol $B_{\delta } \left(0\right)$ is a closed ball of radius $\delta $ centered at the point $0\in R^{n} $. 

 It is clear that the Legendre condition follows from necessary condition \eqref{GrindEQ__1_8_}  as a corollary. We formulate this condition. Let the admissible function $\bar{x}\left(\cdot \right)$ be a weak local minimum in problem \eqref{GrindEQ__1_1_}, \eqref{GrindEQ__1_2_}, and in addition suppose that the integrant $L\left(t,\, x,\, \dot{x}\right)$ is twice differentiable with respect to the variable $\dot{x}$ at the points of the set $\left\{\left(t,\, \bar{x}\left(t\right),\, \dot{\bar{x}}\left(t\right)\right):t\in I\right\}$. Then for all $\xi \in R^{n} $ the following inequalities hold:
\begin{equation} \label{GrindEQ__1_9_} 
\begin{array}{l} {\xi ^{T} \, \bar{L}_{\dot{x}\dot{x}} \left(t\pm \right)\, \xi \ge 0,\, \, \forall t\in \left\{\tau \right\},\, \, \xi ^{T} \, \bar{L}_{\dot{x}\dot{x}} \left(t_{0} +\right)\, \xi \ge 0,\, } \\ {\xi ^{T} \, \bar{L}_{\dot{x}\dot{x}} \, \left(t_{1} -\right)\, \xi \ge 0,\, \, \xi ^{T} \, \bar{L}_{\dot{x}\dot{x}} \left(t\right)\xi \ge 0,\, \, \forall t\in \left(t_{0} ,\, t_{1} \right)\backslash \left\{\tau \right\}.} \end{array} 
\end{equation} 

 The admissible function that satisfies the Euler equation, i.e. equality \eqref{GrindEQ__1_3_} is called an extremal in problem \eqref{GrindEQ__1_1_}, \eqref{GrindEQ__1_2_}. In the classical calculus of variation (see, i.e., \cite{4, 15}), as the main goal, the extremal for a minimum was studied, and various necessary and also sufficient conditions were obtained. Recall that (see [4, 15, 33 and etc.]) in the classical calculus of variations a number of known necessary and also sufficient conditions for a minimum remain powerless in the case when at some point $t=\theta \in I$ for the vector  $\xi=\eta \in R^{n} \backslash \left\{0\right\}$ at least one inequality from \eqref{GrindEQ__1_5_}, \eqref{GrindEQ__1_8_} and \eqref{GrindEQ__1_9_} is fulfilled as an equality, i.e. corresponding necessary condition degenerates at the point  $\theta $ or at the point $\theta $ on the right or left.  Therefore, the study of problem \eqref{GrindEQ__1_1_}, \eqref{GrindEQ__1_2_} in such degenerated cases is of theoretical and practical interest today. 

It should be noted that similar problems in theory of optimal control, starting with the work of Kelly \cite{13}, in terms of singular controls were thoroughly studied by many authors and a hundreds of papers and monographs containing a number of important results were published [2, 7-10, 13, 14, 18-21, 24-28, 32,  and etc.] . Application of the obtained results on singular controls either are ineffective or some their justified modifications are required in degenerated cases when solving problem \eqref{GrindEQ__1_1_}, \eqref{GrindEQ__1_2_}. Although problem  \eqref{GrindEQ__1_1_}, \eqref{GrindEQ__1_2_} is a special case of a terminal optimal control problem with equality type phase constraints, its study as an independent problem, allows to get more effective  results not being corollaries of theorems proved in the theory of optimal control. Confirmation of the last sentence was shown for example in [6, p.33; 22, 23] and in the present paper in the presence of various degenerations. 

In this paper we offer a method for studying an extremum in problem  \eqref{GrindEQ__1_1_}, \eqref{GrindEQ__1_2_} involving degenerations that are based on the introduction of various forms of Weierstrass-type variations characterized by a numerical parameter. Different necessary conditions for a strong and weak local minimum were obtained. It should be noted that some results of this paper are sharpening and refinement of corresponding statements of the paper \cite{3}.  

The structure of this paper is set out as follows. In the second section we obtain increment formulas of the functional for problem \eqref{GrindEQ__1_1_}, \eqref{GrindEQ__1_2_} under various assumptions on the smoothness of the integrant $L\left(\cdot \right)$ and the extremal studied for a minimum (see (2.15), \eqref{GrindEQ__2_16_}, \eqref{GrindEQ__2_32_}, \eqref{GrindEQ__2_33_}, \eqref{GrindEQ__2_36_}, \eqref{GrindEQ__2_37_}).

In the third and fourth sections, based on the increment formula of the  functional obtained in the second section, we introduce necessary conditions for a strong and a weak local minimum in the presence of various degenerations at separate points and on the interval. 

In the last section, by means of special examples we discuss the results obtained in the third and fourth sections.

\section{Various increment formulas of functional in problem \eqref{GrindEQ__1_1_}, \eqref{GrindEQ__1_2_}.}

In this section, by means of special variations we obtain increment formulas of functional under various assumptions on the smoothness of the integrant $L\, \left(\cdot \right)$ and on the considered extremal of problem \eqref{GrindEQ__1_1_}, (1.2). Note that these formulas have independent meaning and are the basis for the proof of the main results of this work.

 Let the admissible function $\bar{x}\, \left(\cdot \right)$ be an extremal in problem \eqref{GrindEQ__1_1_}, \eqref{GrindEQ__1_2_} and  $\vartheta :=\left(\theta ,\, \lambda ,\, \xi \right)\in \left[t_{0} ,\, t_{1} \right)\times \left[0,\, 1\right)\times R^{n} $ be an arbitrary fixed point. Consider a function of the form \cite{22}
\begin{equation} \label{GrindEQ__2_1_} 
h^{\left(+\right)} \left(t;\vartheta ,\, \varepsilon \right)=\left\{\begin{array}{c} {\left(t-\theta \right)\xi ,\, \, } \\ {\frac{\lambda }{\lambda -1} \left(t-\theta -\varepsilon \right)\xi ,} \\ {0,} \end{array}\right. \begin{array}{c} {t\in \left[\theta ,\, \theta +\lambda \varepsilon \right),} \\ {\, \, \, \, \, \, t\in \left[\theta +\lambda \varepsilon ,\, \theta +\varepsilon \right),} \\ {\, t\in I\backslash \left[\theta ,\, \theta +\varepsilon \right).} \end{array} 
\end{equation} 
Here $\varepsilon \in \left(0,\, \varepsilon _{0} \right]$, where $\, \varepsilon _{0} >0$ is a rather small number, moreover $\theta +\varepsilon _{0} <t_{1} $.

Obviously, for any $\varepsilon \in \left(0,\, \varepsilon _{0} \right]$ the function $h^{\left(+\right)\, } \left(\cdot \, ;\, \vartheta ,\, \varepsilon \right)$ is an element of the space $KC^{1} \left(I,R^{n} \right)$ and its derivative $\dot{h}^{\left(+\right)\, } \left(\cdot \, ;\, \vartheta ,\, \varepsilon \right)$ is calculated by 
\begin{equation} \label{GrindEQ__2_2_} 
\dot{h}^{\left(+\right)} \left(t;\, \vartheta ,\, \varepsilon \right)=\left\{\begin{array}{c} {\xi ,\, \, } \\ {\frac{\lambda }{\lambda -1} \xi ,} \\ {0,} \end{array}\right. \begin{array}{c} {t\in \left[\theta ,\, \theta +\lambda \varepsilon \right],} \\ {\, \, \, \, \, \, \, t\in \left[\theta +\lambda \varepsilon ,\, \, \theta +\varepsilon \right],} \\ {\, t\in I\backslash \left(\theta ,\, \theta +\varepsilon \right).} \end{array} 
\end{equation} 

As can be seen, the derivative $\dot{h}^{\left(+\right)\, } \left(\cdot \, ;\, \vartheta ,\, \varepsilon \right)$ at the points $\theta ,\, \theta +\lambda \varepsilon $ and $\theta +\varepsilon $ is calculated both on the right and left, at the points $t_{0} $ and $t_{1} $ on the right and left, respectively. 

Since $\bar{x}\left(\cdot \right)$ is an extremum in problem \eqref{GrindEQ__1_1_}, \eqref{GrindEQ__1_2_}, then by virtue of \eqref{GrindEQ__1_3_}, \eqref{GrindEQ__2_1_} and \eqref{GrindEQ__2_2_} allowing for the property $P\left(\theta +;\, \bar{x}\left(\cdot \right),\, \alpha \right)$ we have the equality   
\begin{equation} \label{GrindEQ__2_3_} 
\begin{array}{c} {0=\int _{\theta }^{\theta +\varepsilon }\frac{d}{dt} \, \left(\bar{L}_{\dot{x}}^{\, T} \left(t\right)\, h_{\varepsilon }^{\left(+\right)} \left(t\right)\right)\, dt= \int _{\theta }^{\theta +\varepsilon }\, \, \, \left[\bar{L}_{x}^{\, T} \left(t\right)\, h_{\varepsilon }^{\left(+\right)} \left(t\right)+\bar{L}_{\dot{x}}^{\, T} \left(t\right)\, \dot{h}_{\varepsilon }^{\left(+\right)} \left(t\right)\right]\, dt ,\, \forall \varepsilon \in \left(0,\, \bar{\varepsilon }\right]} \\ {} \end{array},         
\end{equation} 
where $h_{\varepsilon }^{\left(+\right)} \left(t\right):=h_{\varepsilon }^{\left(+\right)} \left(t;\, \vartheta ,\, \varepsilon \right)$ and $\bar{\varepsilon }=\min \left\{\alpha ,\, \varepsilon _{0} \right\}$.

Note that relation \eqref{GrindEQ__2_3_} is important in the future when calculating the increment of the functional in problem \eqref{GrindEQ__1_1_}, \eqref{GrindEQ__1_2_}. 

Consider the special variation of the function $\bar{x}\left(\cdot \right)$:
\begin{equation} \label{GrindEQ__2_4_} 
x^{\left(+\right)} \left(t;\, \vartheta ,\, \varepsilon \right)=\bar{x}\left(t\right)+h^{\left(+\right)} \left(t;\, \vartheta ,\, \varepsilon \right),\, \, t\in I,\, \, \varepsilon \in \left(0,\, \bar{\varepsilon }\right],   
\end{equation} 
where $h^{\left(+\right)} \left(\cdot \, ;\, \vartheta ,\, \, \varepsilon \right)$ is determined by \eqref{GrindEQ__2_1_}.

 We call \eqref{GrindEQ__2_4_} a variation introduced on the right with respect to the point $\theta $. Obviously, for every $\varepsilon \in \left(0,\, \bar{\varepsilon }\right]$ the function $x^{\left(+\right)} \left(\cdot \, ;\, \vartheta ,\, \, \varepsilon \right)$ is admissible.

Similar to \eqref{GrindEQ__2_4_}, considering the property $P\left(\theta -;\, \bar{x}\left(\cdot \right),\, \alpha \right)$, we introduce into consideration the following variation of the function $\overline{x}\left(\cdot \right)$, the so-called variation introduced on the left with respect to the point $\theta $:
\begin{equation} \label{GrindEQ__2_5_} 
x^{\left(-\right)} \, \left(t;\, \vartheta ,\, \varepsilon \right)=\bar{x}\left(t\right)+h^{\left(-\right)} \left(t:\, \vartheta ,\, \varepsilon \right),\, \, t\in I,\, \, \varepsilon \in \left(0,\, \tilde{\varepsilon }\right],   
\end{equation} 
where $\varepsilon \in \left(0,\, \, \tilde{\varepsilon }\right]$, moreover $\tilde{\varepsilon }=\min \left\{\alpha ,\, \theta -t_{0} \right\}$, further  $\vartheta :=\left(\theta ,\, \lambda ,\, \xi \right)\in \left(t_{0} ,\, t_{1} \, \right]\times \left[0,\, 1\right)\times R^{n} $ is some fixed point and the function $h^{\left(-\right)} \, \left(\cdot \, ;\, \vartheta ,\, \varepsilon \right)$ is determined in the form
\begin{equation} \label{GrindEQ__2_6_} 
h^{\left(-\right)} \, \left(t;\, \vartheta ,\, \varepsilon \right)=\left\{\begin{array}{c} {\left(t-\theta \right)\, \xi ,\, \, } \\ {\frac{\lambda }{\lambda -1} \left(t-\theta +\varepsilon \right)\xi ,} \\ {0,} \end{array}\begin{array}{c} {t\in \left(\theta -\lambda \, \varepsilon ,\, \theta \right]\, ,} \\ {\, \, \, \, \, \, \, t\in \left(\theta -\varepsilon ,\, \, \theta -\lambda \, \varepsilon \right]\, ,} \\ {t\in I\backslash \left(\theta -\varepsilon ,\, \theta \right]\, .} \end{array}\right.  
\end{equation} 

It is clear that for every $\varepsilon \in \left(0,\, \tilde{\varepsilon }\right]$ we have the inclusion $h^{\left(-\right)} \, \left(t;\, \vartheta ,\, \varepsilon \right)\in KC^{1} \left(I,\, R^{n} \right)$ and its derivative is calculated by
\begin{equation} \label{GrindEQ__2_7_} 
\dot{h}^{\left(-\right)} \, \left(t;\, \vartheta ,\, \varepsilon \right)=\left\{\begin{array}{c} {\, \xi ,\, \, } \\ {\frac{\lambda }{\lambda -1} \xi ,} \\ {0,} \end{array}\, \begin{array}{c} {t\in \left[\theta -\lambda \, \varepsilon ,\, \theta \right]\, ,} \\ {\, \, \, \, \, \, t\in \left[\theta -\varepsilon ,\, \, \theta -\lambda \, \varepsilon \right]\, ,} \\ {t\in I\backslash \left(\theta -\varepsilon ,\, \, \theta \right)\, .} \end{array}\right.  
\end{equation} 

By virtue of \eqref{GrindEQ__2_6_} and \eqref{GrindEQ__2_7_} we have that for every $\varepsilon \in \left(0,\, \tilde{\varepsilon }\right]$ the function \eqref{GrindEQ__2_5_} is admissible. 

Similar to \eqref{GrindEQ__2_3_} we confirm that by virtue of \eqref{GrindEQ__1_3_}, \eqref{GrindEQ__2_6_} and \eqref{GrindEQ__2_7_}, allowing for the property $P\left(\theta -;\, \bar{x}\left(\cdot \right),\, \alpha \right)$, the following equality is valid:
\begin{equation} \label{GrindEQ__2_8_} 
0=\int _{\theta -\varepsilon }^{\theta }\frac{d}{dt} \, \left(\bar{L}_{\dot{x}}^{\, T} \left(t\right)\, h_{\varepsilon }^{\left(-\right)} \left(t\right)\right)\, dt= \int _{\theta -\varepsilon }^{\theta }\left[\bar{L}_{x}^{\, T} \left(t\right)\, h_{\varepsilon }^{-} \left(t\right)+\bar{L}_{\dot{x}}^{\, T} \left(t\right)\, \dot{h}_{\varepsilon }^{\left(-\right)} \left(t\right)\right]dt, \, \forall \varepsilon \in \left(0,\, \tilde{\varepsilon }\right],  
\end{equation} 
where $h_{\varepsilon }^{\left(-\right)} \left(t\right):=h^{\left(-\right)} \left(t;\, \vartheta ,\, \varepsilon \right)$.

 We introduce denotations that will be convenient in what follows:
\begin{equation} \label{GrindEQ__2_9_} 
\bar{L}_{x} \left(t,\, \xi \right):=L_{x} \left(t,\bar{x}\left(t\right),\, \dot{\bar{x}}\left(t\right)+\xi \right),\, \, \, \, \bar{L}_{xx} \left(t,\, \xi \right):=L_{xx} \left(t,\bar{x}\left(t\right),\, \dot{\bar{x}}\left(t\right)+\xi \right)\, ,   
\end{equation} 
\begin{equation} \label{GrindEQ__2_10_} 
\Delta \bar{L}_{x} \left(t,\, \xi \right):=L_{x} \left(t,\bar{x}\left(t\right),\, \dot{\bar{x}}\left(t\right)+\xi \right)-\, \bar{L}_{x} \left(t\right),\, \xi \in R^{n} ,      
\end{equation} 
\begin{equation} \label{GrindEQ__2_11_} 
Q_{i} \left(\bar{L}\right)\, \left(t,\, \lambda \, ,\, \xi \right):=\lambda ^{i} E\, \left(\bar{L}\right)\, \left(t,\, \, \xi \right)+\left(1-\lambda ^{i} \right)E\left(\bar{L}\right)\, \left(t,\, \frac{\lambda }{\lambda -1} \xi \right),\, \, i=1,2,3,   
\end{equation} 
\begin{equation} \label{GrindEQ__2_12_} 
M_{i} \left(\bar{L}_{x} \right)\left(t,\, \lambda ,\, \xi \right):=\lambda ^{i} \Delta \, \bar{L}_{x}^{\, T} \left(t,\, \, \xi \right)\xi +\left(1-\lambda \right)\left(\frac{1}{2} +\lambda \right)^{i-1} \Delta \, \bar{L}_{x}^{\, T} \left(t,\, \frac{\lambda }{\lambda -1} \xi \right)\xi ,\, \, \, i=1,2. 
\end{equation} 

Considering the above stated, we prove the following propositions. 

\textbf{Proposition 2.1. }Let the functions $L\left(\cdot \right)$ and $L_{\dot{x}} \left(\cdot \right)$ be twice continuously differentiable in totality of variables, and the admissible function $\bar{x}\left(\cdot \right)$ be an extremal in problem \eqref{GrindEQ__1_1_}, \eqref{GrindEQ__1_2_}. Then:

(i) if at the point $\theta \in \left[t{}_{0} ,\, t{}_{1} \, \right)$ the extremal $\bar{x}\left(\cdot \right)$ is triply differentiable on the right in right semi-neighborhood $\left[\theta ,\, \theta +\alpha \right)\subset I$ of the point $\theta $, then for every $\left(\lambda ,\, \xi \right)\in \, \left[0,1\right)\times R^{n} $ there exists such a number $\varepsilon ^{*} >0$ that for all $\varepsilon \in \left(0,\, \varepsilon ^{*} \right]$ the increment $J\left(x^{\left(+\right)} \left(\cdot :\, \vartheta ,\varepsilon \right)\right)-J\left(\bar{x}\left(\cdot \right)\right)=:\Delta _{\varepsilon }^{\left(+\right)} J\left(\bar{x}\left(\cdot \right);\vartheta \right)$ of the functional \eqref{GrindEQ__1_1_}, corresponding to the variation \eqref{GrindEQ__2_4_} is represented in the form 
\[\Delta _{\varepsilon }^{\left(+\right)} J\left(\bar{x}\left(\cdot \right);\, \vartheta \right)\, =\varepsilon \, Q_{1} \left(\bar{L}\right)\, \left(\theta +,\, \lambda ,\, \xi \right)+\frac{1}{2} \varepsilon ^{2} W\left(\bar{L}\right)\, \left(\theta +,\, \lambda ,\, \xi \right)+\] 
\begin{equation} \label{GrindEQ__2_13_} 
+\frac{1}{6} \varepsilon ^{3} G\left(\bar{L}\right)\, \left(\theta +,\, \lambda ,\, \xi \right)+o\left(\varepsilon ^{3} \right);    
\end{equation} 

(ii) if at the point $\theta \in \left(t{}_{0} ,\, t{}_{1} \right]$ the extremal $\bar{x}\left(\cdot \right)$ is triply differentiable on the left in left semi-neighborhood $\left(\theta -\alpha ,\, \theta \right]\subset I$ of the point $\theta $, then for every $\left(\lambda ,\, \xi \right)\in \left[0,\, 1\right)\times R^{n} $ there exists such a number $\varepsilon ^{*} >0$ that for all $\varepsilon \in \left(0,\, \varepsilon ^{*} \right]$ the increment $J\left(x^{\left(-\right)} \left(\cdot ,\, \vartheta ,\, \varepsilon \right)\right)-J\left(\bar{x}\left(\cdot \right)\right)=:\Delta _{\varepsilon }^{\left(-\right)} J\left(\bar{x}\left(\cdot \right);\, \vartheta \right)$ corresponding to the variation \eqref{GrindEQ__2_5_}, is represented in the form
\[\Delta _{\varepsilon }^{\left(-\right)} J\left(\bar{x}\left(\cdot \right);\, \vartheta \right)\, =\varepsilon \, Q_{1} \left(\bar{L}\right)\, \left(\theta -,\, \lambda ,\, \xi \right)-\frac{1}{2} \varepsilon ^{2} W\left(\bar{L}\right)\, \left(\theta -,\, \lambda ,\, \xi \right)+\] 
\begin{equation} \label{GrindEQ__2_14_} 
+\frac{1}{6} \varepsilon ^{3} G\left(\bar{L}\right)\, \left(\theta -,\, \lambda ,\, \xi \right)+o\left(\varepsilon ^{3} \right).    
\end{equation} 
Here
\begin{equation} \label{GrindEQ__2_15_} 
W\left(\bar{L}\right)\, \left(\theta ,\, \lambda ,\, \xi \right):=\lambda \, M_{1} \left(\bar{L}_{x} \right)\, \left(\theta ,\, \lambda ,\, \xi \right)+\frac{d}{dt} Q_{2} \left(\bar{L}\right)\, \left(\theta ,\, \lambda ,\, \xi \right),   
\end{equation} 
\[G\left(\bar{L}\right)\, \left(\theta ,\, \lambda ,\, \xi \right):=\lambda ^{2} \xi ^{T} \left[\lambda \bar{L}_{xx} \left(\theta ,\, \xi \right)+\left(1-\lambda \right)\bar{L}_{xx} \left(\theta ,\, \frac{\lambda }{\lambda -1} \xi \right)\right]\xi +\] 
\begin{equation} \label{GrindEQ__2_16_} 
+2\lambda \frac{d}{dt} M_{2} \left(\bar{L}_{x} \right)\, \left(\theta ,\, \lambda ,\, \xi \right)+\frac{d^{2} }{dt^{2} } Q_{3} \left(\bar{L}\right)\, \left(\theta ,\, \lambda ,\, \xi \right),    
\end{equation} 
where $Q_{i} \left(\bar{L}\right)\, \left(\cdot \right),\, \, i=1,2,3$ and $M_{i} \left(\bar{L}_{x} \right)\, \left(\cdot \right),\, \, i=1,2$  are determined by \eqref{GrindEQ__2_11_} and \eqref{GrindEQ__2_12_}, allowing for \eqref{GrindEQ__1_6_}, \eqref{GrindEQ__2_9_} and \eqref{GrindEQ__2_10_}.

 \textbf{Proof. }At first we prove part (i) of the proposition, i.e. the validity of the equality \eqref{GrindEQ__2_13_}.

Since the admissible function $\bar{x}\left(\cdot \right)$ is an extremal in problem \eqref{GrindEQ__1_1_}, \eqref{GrindEQ__1_2_}, i.e. is the solution of the equation\eqref{GrindEQ__1_3_}, then using  \eqref{GrindEQ__2_3_} and \eqref{GrindEQ__2_4_}, allowing for \eqref{GrindEQ__2_1_}, \eqref{GrindEQ__2_2_} and property $P\left(\theta +,\, \bar{x}\left(\cdot \right),\, \alpha \right)$, it is easy to represent the increment $\Delta _{\varepsilon }^{\left(+\right)} \, J\left(\bar{x}\, \left(\cdot \right),\, \vartheta \right)$ in the form:
\begin{equation} \label{GrindEQ__2_17_} 
\Delta _{\varepsilon }^{\left(+\right)} \, J\left(\bar{x}\, \left(\cdot \right),\, \vartheta \right)=J_{1}^{\left(+\right)} \left(\varepsilon ,\, \vartheta \right)+J_{2}^{\left(+\right)} \left(\varepsilon ,\, \vartheta \right),\, \, \varepsilon \in \left(0,\, \bar{\varepsilon }\right].        
\end{equation} 
Here 
\begin{equation} \label{GrindEQ__2_18_} 
\begin{split}
J_{1}^{\left(+\right)} \left(\varepsilon ,\, \vartheta \right)=\int _{\theta }^{\theta +\varepsilon }\, \, \, \left[L\left(t,\, \bar{x}\left(t\right),\, \dot{\bar{x}}\left(t\right)+\dot{h}{}_{\varepsilon }^{\left(+\right)} \left(t\right)\right)-\bar{L}\left(t\right)-\bar{L}_{\dot{x}}^{T} \left(t\right)\dot{h}{}_{\varepsilon }^{\left(+\right)} \left(t\right)\right] \, \, dt\, ,  
\end{split} 
\end{equation} 

\begin{equation} \label{GrindEQ__2_19_} 
\begin{split}
J_{2}^{\left(+\right)} \left(\varepsilon ,\, \vartheta \right)=\int _{\theta }^{\theta +\varepsilon }\, \, \, [L(t,\, \bar{x}\left(t\right)+h{}_{\varepsilon }^{\left(+\right)} \left(t\right),\dot{\bar{x}}\left(t\right)+\dot{h}{}_{\varepsilon }^{\left(+\right)} (t))\\
- L\, \left(t,\, \bar{x}\left(t\right),\, \dot{\bar{x}}\left(t\right)+\dot{h}{}_{\varepsilon }^{\left(+\right)} \left(t\right)\right)-\bar{L}_{x}^{T} \left(t\right)h{}_{\varepsilon }^{\left(+\right)} \left(t\right)]dt,   
\end{split}
\end{equation} 
where $h{}_{\varepsilon }^{\left(+\right)} \left(t\right):=h^{\left(+\right)} \left(t;\, \vartheta ,\, \varepsilon \right)$, and the number $\bar{\varepsilon }$ is determined above (see \eqref{GrindEQ__2_3_}).

 We calculate the integrals \eqref{GrindEQ__2_18_} and \eqref{GrindEQ__2_19_} with accuracy $o(\varepsilon ^{3})$. Considering the assumption on the smoothness of the functions $L\, \, \left(\cdot \right)$ and $\bar{x}\left(\cdot \right)$, we apply the Taylor formula. Then:

(a) by virtue of \eqref{GrindEQ__2_1_} - \eqref{GrindEQ__2_3_}, allowing for  \eqref{GrindEQ__1_6_}, \eqref{GrindEQ__1_7_} and \eqref{GrindEQ__2_11_}, from \eqref{GrindEQ__2_18_} we have 
\[J_{1}^{\left(+\right)} \left(\varepsilon ,\, \vartheta \right)=\int _{\theta }^{\theta +\lambda \varepsilon }\, E\left(\bar{L}\right)\, \left(t,\, \xi \right) \, dt\, +\int _{\theta +\lambda \varepsilon }^{\theta +\varepsilon }E\left(\bar{L}\right)\, \left(t,\, \frac{\lambda }{\lambda -1} \xi \right)dt= \] 
\[=\int _{\theta }^{\theta +\lambda \varepsilon }\left[\sum _{i=0}^{2}\frac{1}{i!}  \left(\, t-\theta \right)^{i} \frac{d^{i} }{dt^{i} } E\left(\bar{L}\right)\, \left(\theta +,\xi \right)\, +o\left(\left(t-\theta \right)^{2} \right)\right]\, dt+ \] 
\[+\int _{\theta +\lambda \varepsilon }^{\theta +\varepsilon }\left[\sum _{i=0}^{2}\frac{1}{i!}  \left(t-\theta  \right)^{i}\frac{d^{i} }{dt^{i} } E\left(\bar{L}\right)\, \left(\theta +,\, \frac{\lambda }{\lambda -1} \xi \right)+o\left(\left(t-\theta \right)^{2} \right)\right]dt= \] 
\[=\sum _{i=0}^{2}\frac{\varepsilon ^{i+1} }{\left(i+1\right)\, !}  \frac{d^{i} }{dt^{i} } \left[\lambda ^{i+1} E\left(\bar{L}\right)\, \left(\theta +,\xi \right)+\left(1-\lambda ^{i+1} \right)E\left(\bar{L}\right)\, \left(\theta +,\, \frac{\lambda }{\lambda -1} \xi \right)\right]+\] 
\begin{equation} \label{GrindEQ__2_20_} 
+o\left(\varepsilon ^{3} \right)=\, \sum _{i=0}^{2}\frac{\varepsilon ^{i+1} }{\left(i+1\right)!}  \frac{d^{i} }{dt^{i} } Q_{i+1} \left(\bar{L}\right)\, \left(\theta +,\, \lambda ,\, \xi \right)+o\left(\varepsilon ^{3} \right);   
\end{equation} 

(b) by virtue of \eqref{GrindEQ__2_1_} - \eqref{GrindEQ__2_3_}, allowing for \eqref{GrindEQ__2_9_} and \eqref{GrindEQ__2_10_}, from \eqref{GrindEQ__2_19_} we get 
\begin{equation} \label{GrindEQ__2_21_} 
J{}_{2}^{\left(+\right)} \left(\varepsilon ,\, \vartheta \right)=J{}_{21}^{\left(+\right)} \left(\varepsilon ,\, \vartheta \right)+J{}_{22}^{\left(+\right)} \left(\varepsilon ,\, \vartheta \right),\, \, \forall \varepsilon \in \left(0,\, \bar{\varepsilon }\right].   
\end{equation} 
Here the integrals $J{}_{21}^{\left(+\right)} \left(\cdot \right)$ and $J{}_{22}^{\left(+\right)} \left(\cdot \right)$ are calculated by the Taylor formula as follows: 
\[J{}_{21}^{\left(+\right)} \left(\varepsilon ,\, \vartheta \right)=\int _{\theta }^{\theta +\lambda \varepsilon }\left[\left(t-\theta \right)\Delta \, \bar{L}_{x}^{\, T} \left(t,\, \xi \right)\xi +\frac{1}{2} \left(t-\theta \right)^{2} \xi ^{T} \bar{L}_{xx} \left(t,\, \xi \right)\xi +o\left(\left(t-\theta \right)^{2} \right)\, \right]\, dt =\] 
\[=\int _{\theta }^{\theta +\lambda \varepsilon }\left(t-\theta \right)\, \left[\Delta \bar{L}_{x}^{\, T} \left(\theta +,\, \xi \right)\, \xi +\left(t-\theta \right)\frac{d}{dt} \Delta \, \bar{L}_{x}^{\, T} \left(\theta +,\, \xi \right)\, \xi +o\left(t-\theta \right)\right] \, dt+\] 
\[+\frac{1}{6} \varepsilon ^{3} \lambda ^{3} \xi ^{T} \bar{L}_{xx} \left(\theta +,\, \xi \right)\xi +o\, \left(\varepsilon ^{3} \right)=\frac{\varepsilon ^{2} }{2} \lambda ^{2} \Delta \, \bar{L}_{x}^{\, T} \left(\theta +,\, \xi \right)\, \xi +\] 
\begin{equation} \label{GrindEQ__2_22_} 
+\frac{\varepsilon ^{3} }{6} \, \left[2\lambda ^{3} \frac{d}{dt} \Delta \, \bar{L}_{x}^{\, T} \left(\theta +,\, \xi \right)\xi +\lambda ^{3} \xi ^{T} \bar{L}_{xx} \left(\theta +,\, \xi \right)\xi \right]+o\left(\varepsilon ^{3} \right),   
\end{equation} 
\[J{}_{22}^{\left(+\right)} \left(\varepsilon ,\, \vartheta \right)=\int _{\theta +\lambda \varepsilon }^{\theta +\varepsilon }q\left(t;\theta ,\, \lambda ,\, \varepsilon \right)dt ,     \] 
where 
\[q\left(t;\theta ,\, \lambda ,\, \varepsilon \right)=\left(\lambda -1\right)^{-1} \lambda \, \left(t-\theta -\varepsilon \right)\Delta \bar{L}_{x}^{T} \left(t,\, \left(\lambda -1\right)^{-1} \lambda \xi \right)\xi +\] 
\[+\frac{1}{2} \left(\lambda -1\right)^{-2} \lambda ^{2} \left(t-\theta -\varepsilon \right)^{2} \xi ^{T} \bar{L}_{xx} \left(t,\, \left(\lambda -1\right)^{-1} \lambda \, \xi \right)\, \xi +o\, \left(\left(t-\theta -\varepsilon \right)^{2} \right).      \] 
Continuing the calculations, we have 
\[J{}_{22}^{\left(+\right)} \left(\varepsilon ,\, \vartheta \right)=\] 
\[=\frac{\lambda }{\lambda -1} \int _{\theta +\lambda \varepsilon }^{\theta +\varepsilon }\left(t-\theta -\varepsilon \right)\left[\Delta \bar{L}_{x} \left(\theta +,\, \frac{\lambda }{\lambda -1} \xi \right)\xi +\left(t-\theta \right)\frac{d}{dt} \Delta \bar{L}_{x} \left(\theta +,\, \frac{\lambda }{\lambda -1} \xi \right)\xi +o\left(t-\theta \right)\right] dt+\] 
\[+\frac{1}{6} \varepsilon ^{3} \lambda ^{2} \left(1-\lambda \right)\xi ^{T} \bar{L}_{xx} \left(\theta +,\, \frac{\lambda }{\lambda -1} \xi \right)\xi +o \left(\varepsilon ^{3} \right)=\frac{\varepsilon ^{2} }{2} \lambda \left(1-\lambda \right)\Delta \bar{L}_{x}^{T} \left(\theta +,\, \frac{\lambda }{\lambda -1} \xi \right)\xi +\] 
\[+\frac{\varepsilon ^{3} }{6} \left[\lambda \left(1-\lambda \right)\left(1+2\lambda \right)\frac{d}{dt} \Delta \bar{L}_{x}^{T} \left(\theta +,\, \frac{\lambda }{\lambda -1} \xi \right)\xi +\left. \lambda ^{2} \left(1-\lambda \right)\xi ^{T} \bar{L}_{xx} \left(\theta +,\, \frac{\lambda }{\lambda -1} \xi \right)\xi \right]+o\left(\varepsilon ^{3} \right).\right . \] 
\begin{equation} \label{GrindEQ__2_23_} 
.    
\end{equation} 

By virtue of \eqref{GrindEQ__2_22_} and \eqref{GrindEQ__2_23_}, allowing for the notation \eqref{GrindEQ__2_12_}, the equality \eqref{GrindEQ__2_21_}, takes the form
\[J_{2}^{\left(+\right)} \left(\varepsilon ,\, \vartheta \right)=\frac{1}{2} \varepsilon ^{2} \lambda \, M_{1} \left(\bar{L}_{x} \right)\, \left(\theta +,\, \lambda ,\xi \right)+\frac{1}{6} \varepsilon ^{3} \, \left[\lambda ^{2} \xi ^{T} \right. \lambda \, \bar{L}_{xx} \left(\theta +,\, \xi \right)+\] 
\begin{equation} \label{GrindEQ__2_24_} 
\left. +\left(1-\lambda \right)\bar{L}_{xx} \left(\theta +,\left(\lambda -1\right)^{-1} \lambda \xi \right)\xi +2\lambda \frac{d}{dt} M_{2} \left(\bar{L}_{x} \right)\, \left(\theta +,\, \lambda ,\xi \right)\right]+o\left(\varepsilon ^{3} \right).  
\end{equation} 

Consequently, having substituted \eqref{GrindEQ__2_20_} and \eqref{GrindEQ__2_24_} in \eqref{GrindEQ__2_17_}, and also having chosen $\varepsilon ^{*} =\bar{\varepsilon }$, allowing for \eqref{GrindEQ__2_10_}-\eqref{GrindEQ__2_12_}, \eqref{GrindEQ__2_15_} and \eqref{GrindEQ__2_16_}, we get the expansion \eqref{GrindEQ__2_13_}, i.e. part (i) of proposition 2.1 is proved.  

Similar to  \eqref{GrindEQ__2_13_} we give the proof of part  (ii) of proposition 2.1. 

For that it is sufficient to calculate the increment  $J\left(x^{\left(-\right)} \left(\cdot ;\, \vartheta ,\varepsilon \right)\right)-J\left(\bar{x}\left(\cdot \right)\right)=:\Delta _{\varepsilon }^{\left(-\right)} J\left(\bar{x}\left(\cdot \right),\, \vartheta \right)$ with accuracy $o(\varepsilon ^{3} )$, where the function $x^{(-)} (\cdot ;\, \vartheta ,\, \varepsilon )$ is determined by \eqref{GrindEQ__2_5_} and \eqref{GrindEQ__2_6_}. Using \eqref{GrindEQ__2_5_} and \eqref{GrindEQ__2_8_}, allowing for \eqref{GrindEQ__2_6_}, \eqref{GrindEQ__2_7_} and the property $P\left(\theta -,\, \bar{x}\left(\cdot \right),\, \alpha \right)$, we can represent the increment $\Delta _{\varepsilon }^{\left(-\right)} J\left(\bar{x}\, \left(\cdot \right),\, \vartheta \right)$ in the form 
\begin{equation} \label{GrindEQ__2_25_} 
\Delta _{\varepsilon }^{\left(-\right)} J\, \left(\bar{x}\, \left(\cdot \right),\, \vartheta \right)=J_{1}^{\left(-\right)} \left(\varepsilon ,\vartheta \right)+J_{2}^{\left(-\right)} \left(\varepsilon ,\vartheta \right),\, \, \varepsilon \in \left(0,\, \tilde{\varepsilon }\right].   
\end{equation} 
Here the number $\tilde{\varepsilon }>0$ is determined above (see \eqref{GrindEQ__2_5_}).
\begin{equation} \label{GrindEQ__2_26_} 
J_{1}^{\left(-\right)} \left(\varepsilon ,\, \vartheta \right):=\int _{\theta -\varepsilon }^{\theta }\, \left[L\, \left(t,\, \bar{x}\left(t\right),\, \dot{\bar{x}}\left(t\right)+\dot{h}_{\varepsilon }^{\left(-\right)} \left(t\right)\right)-\bar{L}\left(t\right)-\bar{L}_{\dot{x}}^{T} \left(t\right)\, \dot{h}_{\varepsilon }^{\left(-\right)} \left(t\right)\right]\, dt ,  
\end{equation} 
\[J_{2}^{\left(-\right)} \left(\varepsilon ,\, \vartheta \right):=\] 
\begin{equation} \label{GrindEQ__2_27_} 
\int _{\theta -\varepsilon }^{\theta }\, \left[L\, \left(t,\, \bar{x}\left(t\right)+h_{\varepsilon }^{\left(-\right)} \left(t\right),\, \, \dot{\bar{x}}\left(t\right)+\dot{h}_{\varepsilon }^{\left(-\right)} \left(t\right)\right)-\right.  \left. L\, \left(t,\, \bar{x}\left(t\right),\, \dot{\bar{x}}\left(t\right)+\dot{h}_{\varepsilon }^{\left(-\right)} \left(t\right)\right)-\bar{L}_{x}^{T} \left(t\right)\, h_{\varepsilon }^{\left(-\right)} \left(t\right)\right]\, dt,   
\end{equation} 
where $h_{\varepsilon }^{\left(-\right)} \left(t\right):=h^{\left(-\right)} \left(t;\, \vartheta ,\, \varepsilon \right)$.

Considering the assumption on the smoothness of functions $L\left(\cdot \right)$,$L_{\dot{x}} \, \left(\cdot \right)$ and $\bar{x}\left(\cdot \right)$, having applied the Taylor formula, we calculate the integrals \eqref{GrindEQ__2_26_} and \eqref{GrindEQ__2_27_} with accuracy $o(\varepsilon ^{3})$. More exactly, we carry out calculations in the following way:  

(a) similar to \eqref{GrindEQ__2_20_}, by virtue of \eqref{GrindEQ__2_6_} and \eqref{GrindEQ__2_7_}, allowing for notations  \eqref{GrindEQ__1_6_} and \eqref{GrindEQ__2_11_}, from \eqref{GrindEQ__2_26_} by the Taylor formula we get

\begin{equation} \label{GrindEQ__2_28_} 
\begin{split}
J_{1}^{\left(-\right)} \left(\varepsilon ,\, \vartheta \right)=\int _{\theta -\lambda \varepsilon }^{\theta }\left[\sum_{i=0}^{2} \frac{1}{i!}(t-\theta)^i \frac{d^i}{dt^i}E(\overline{L})(\theta-, \xi)+o((t-\theta)^2)\right]dt\\
+\int _{\theta -\varepsilon}^{\theta -\lambda \varepsilon }\left[\sum_{i=0}^{2} \frac{1}{i!}(t-\theta)^i \frac{d^i}{dt^i}E(\overline{L})(\theta-, \frac{\lambda}{\lambda-1}\xi)+o((t-\theta)^2)\right]dt \\
=\sum _{i=0}^{2}\left(-1\right)^{i+2} \frac{\varepsilon ^{i+1} }{\left(i+1\right)!} \frac{d^{i} }{dt^{i} } Q_{i+1} \left(\bar{L}\right)\, \left(\theta -,\, \lambda ,\, \xi \right)+o\left(\varepsilon ^{{\rm 3}} \right) ,\, \, \forall \varepsilon \in \left(0,\, \tilde{\varepsilon }\right];   
\end{split}
\end{equation} 

(b) similar to \eqref{GrindEQ__2_22_} and \eqref{GrindEQ__2_23_}, by virtue of \eqref{GrindEQ__2_6_} and \eqref{GrindEQ__2_7_}, allowing for notations \eqref{GrindEQ__2_9_} and \eqref{GrindEQ__2_10_}, from \eqref{GrindEQ__2_27_}  by the Taylor formula we have
\begin{equation} \label{GrindEQ__2_29_} 
J_{2}^{\left(-\right)} \left(\varepsilon ,\, \vartheta \right)=J_{21}^{\left(-\right)} \left(\varepsilon ,\, \vartheta \right)+J_{22}^{\left(-\right)} \left(\varepsilon ,\, \vartheta \right),\, \, \varepsilon \in \left(0,\, \tilde{\varepsilon }\right],    
\end{equation} 
where
\[J_{21}^{\left(-\right)} \left(\varepsilon ,\, \vartheta \right)=\int _{\theta -\lambda \varepsilon }^{\theta }\left[\left(t-\theta \right)\Delta \bar{L}_{x}^{T} \, \left(t,\, \xi \right)\, \xi +\frac{1}{2} \left(t-\theta \right)^{2} \xi ^{T} \bar{L}_{xx} \left(t,\, \xi \right)\xi +o\left(\left(t-\theta \right)^{2} \right)\right]dt= \] 
\begin{equation} \label{GrindEQ__2_30_} 
=-\frac{\varepsilon ^{2} }{2} \lambda ^{2} \Delta \bar{L}_{x}^{T} \, \left(\theta -,\, \xi \right)\, \xi +\frac{1}{6} \varepsilon ^{3} \left. \left[2\lambda ^{3} \frac{d}{dt} \Delta \bar{L}_{x}^{T} \, \left(\theta -,\, \xi \right)\, \xi \right. +\lambda ^{3} \xi ^{T} \bar{L}_{xx} \left(\theta -,\, \xi \right)\, \xi \right]+o\left(\varepsilon ^{3} \right), 
\end{equation} 
\[J_{22}^{\left(-\right)} \left(\varepsilon ,\, \vartheta \right)=\int _{\theta -\varepsilon }^{\theta -\lambda \varepsilon }\frac{\lambda }{\lambda -1} \left(t-\theta +\varepsilon \right)\Delta \bar{L}_{x}^{T} \, \left(t,\, \frac{\lambda }{\lambda -1} \xi \right)\xi \, dt+ \] 
\[+\int _{\theta -\varepsilon }^{\theta -\lambda \varepsilon }\left[\frac{\lambda ^{2} }{2\left(\lambda -1\right)^{2} } \left(t-\theta +\varepsilon \right)^{2} \xi ^{T} \bar{L}_{xx} \left(t,\, \frac{\lambda }{\lambda -1} \xi \right)\xi +o\left(\left(t-\theta +\varepsilon \right)^{2} \right)\right]\, dt= \] 
\[=-\frac{1}{2} \varepsilon ^{2} \lambda \left(1-\lambda \right)\Delta \, \bar{L}_{x}^{T} \, \left(\theta -,\, \frac{\lambda }{\lambda -1} \xi \right)\xi +\] 
\[+\frac{\varepsilon ^{3} }{6} \left[\lambda \left(1-\lambda \right)\, \left(1+2\lambda \right)\frac{d}{dt} \Delta \, \bar{L}_{x}^{T} \, \left(\theta -,\, \frac{\lambda }{\lambda -1} \xi \right)\xi +\lambda ^{2} \left(1-\lambda \right)\xi ^{T} \bar{L}_{xx} \left(\theta -,\, \frac{\lambda }{\lambda -1} \xi \right)\xi \right]+o\left(\varepsilon ^{3} \right).\]

Consequently, by virtue of \eqref{GrindEQ__2_28_}-\eqref{GrindEQ__2_30_} for the increment \eqref{GrindEQ__2_25_}, taking into account \eqref{GrindEQ__2_10_}-\eqref{GrindEQ__2_12_}, \eqref{GrindEQ__2_15_} and \eqref{GrindEQ__2_16_}, and also choosing $\varepsilon ^{*} =\tilde{\varepsilon }$, we get expansion \eqref{GrindEQ__2_14_}, i.e. part (ii) of Proposition 2.1 is proved. By the same token Proposition 2.1 is completely proved.

Based on the technique for proving Proposition 2.1, namely, using \eqref{GrindEQ__2_17_}-\eqref{GrindEQ__2_19_} and \eqref{GrindEQ__2_25_}-\eqref{GrindEQ__2_27_}, under weak assumptions on the smoothness of functions $L\left(\cdot \right),\, L_{\dot{x}} \left(\cdot \right)$ and $\bar{x}\, \, \left(\cdot \right)$ as a corollary of Proposition 2.1, it is easy to arrive at the following statement.

\textbf{Proposition 2.2. }Let the functions $L\left(\cdot \right)$ and $L_{\dot{x}} \left(\cdot \right)$ be continuously differentiable in totality of variables, and the admissible function $\bar{x}\, \, \left(\cdot \right)$ be an extremal in problem \eqref{GrindEQ__1_1_}, \eqref{GrindEQ__1_2_}. If at the point $\theta \in \left[t_{0} ,\, t_{1} \right)$ $\left(\theta \in \left(t_{0} ,\, t_{1} \right]\right)$ the extremal $\bar{x}\, \, \left(\cdot \right)$ is twice differentiable on the right (left) in semi-neighborhood $\left[\theta ,\, \theta +\alpha \right)\subset I$$\left(\, \left(\theta -\alpha ,\, \theta \right]\subset I\, \right)$ of the point $\theta $, then for every $\left(\lambda ,\, \xi \right)\in \left[0,\, 1\right)\times R^{n} $ there exists such a number $\varepsilon ^{*} >0$ that for all $\varepsilon \in \left(0,\, \varepsilon ^{*} \right]$ the increment $\Delta _{\varepsilon }^{\left(+\right)} \, J\left(\bar{x}\, \left(\cdot \right);\vartheta \right)$ $\left(\Delta _{\varepsilon }^{\left(-\right)} \, J\left(\bar{x}\, \left(\cdot \right);\vartheta \right)\right)$ of the functional in problem \eqref{GrindEQ__1_1_}, \eqref{GrindEQ__1_2_}, corresponding to the variation \eqref{GrindEQ__2_4_} ((2.5)), is represented in the form  
\begin{equation} \label{GrindEQ__2_32_} 
\Delta _{\varepsilon }^{\left(+\right)} J\left(\bar{x}\left(\cdot \right); \vartheta \right)\, :=\varepsilon \, Q_{1} \left(\bar{L}\right)\, \left(\theta +,\, \lambda ,\, \xi \right)+\frac{1}{2} \varepsilon ^{2} W\left(\bar{L}\right)\, \left(\theta +,\, \lambda ,\, \xi \right)+o\left(\varepsilon ^{2} \right) 
\end{equation} 
\begin{equation} \label{GrindEQ__2_33_} 
\left(\Delta _{\varepsilon }^{\left(-\right)} J\left(\bar{x}\left(\cdot \right);\, \vartheta \right)\, :=\varepsilon \, Q_{1} \left(\bar{L}\right)\, \left(\theta -,\, \lambda ,\, \xi \right)-\frac{1}{2} \varepsilon ^{2} W\left(\bar{L}\right)\, \left(\theta -,\, \lambda ,\, \xi \right)+o\left(\varepsilon ^{2} \right)\right),   
\end{equation} 
where $Q_{1}(\bar{L})(\cdot)$ and $W(\bar{L})(\cdot)$ are determined by \eqref{GrindEQ__2_11_} and \eqref{GrindEQ__2_15_} allowing for \eqref{GrindEQ__1_6_} and \eqref{GrindEQ__2_9_}-\eqref{GrindEQ__2_12_}.

 We now consider the following special case. Namely, assuming  $\lambda =\varepsilon \in \left(0,1\right)\bigcap \left(0,\, \bar{\varepsilon }\right]$ $\left(\lambda =\varepsilon \in \left(0,1\right)\bigcap \left(0,\, \tilde{\varepsilon }\right]\right)$ in \eqref{GrindEQ__2_1_} ((2.6)), we introduce a new variation of the extremal $\bar{x}(\cdot)$ in the form: 
\begin{equation} \label{GrindEQ__2_34_} 
\left. x^{\left(+\right)} \left(t;\, \tilde{\vartheta },\, \varepsilon \right):=x^{\left(+\right)} \left(t;\, \vartheta ,\, \varepsilon \right)\, \right|_{\vartheta =\tilde{\vartheta }=\left(\theta ,\varepsilon ,\xi \right)} ,\, \varepsilon \in \left(0,1\right)\bigcap \left(0,\bar{\varepsilon }\right] 
\end{equation} 
\begin{equation} \label{GrindEQ__2_35_} 
\left(\left. x^{\left(-\right)} \left(t;\, \tilde{\vartheta },\, \varepsilon \right):=x^{\left(-\right)} \left(t;\, \vartheta ,\, \varepsilon \right)\, \right|_{\vartheta =\tilde{\vartheta }=\left(\theta ,\varepsilon ,\xi \right)} ,\, \varepsilon \in \left(0,1\right)\bigcap \left(0,\tilde{\varepsilon }\right]\right),   
\end{equation} 
where $\vartheta =\left(\theta ,\, \lambda ,\, \xi \right)$ and $x^{\left(+\right)} \left(\cdot ;\, \vartheta ,\, \varepsilon \right)$ $\left(x^{\left(-\right)} \left(\cdot ;\, \vartheta ,\, \varepsilon \right)\right)$ is determined by \eqref{GrindEQ__2_4_} allowing for \eqref{GrindEQ__2_1_}      (by (2.5) allowing for \eqref{GrindEQ__2_6_}). In this case the following proposition is valid. 

 \textbf{Proposition 2.3.  }Let the functions $L\left(\cdot \right)$ and $L_{\dot{x}} \left(\cdot \right)$ be twice differentiable in totality of variables, the function $\bar{x}\, \, \left(\cdot \right)$ be an extremal in problem \eqref{GrindEQ__1_1_}, \eqref{GrindEQ__1_2_}. If at the point $\theta \in \left[t_{0} ,\, t_{1} \right)$ $\left(\theta \in \left(t_{0} ,\, t_{1} \right]\right)$ the extremal $\bar{x}\, \, \left(\cdot \right)$ is twice differentiable on the right (left) in some semi-neighborhood $\left[\theta ,\, \theta +\alpha \right)\subset I$$\left(\, \left(\theta -\alpha ,\, \theta \right]\subset I\, \right)$ of the point $\theta $, then for every $\xi \in R^{n} $ there exists such a number $\varepsilon ^{*} >0$ that for all $\varepsilon \in \left(0,\, \varepsilon ^{*} \right]\bigcap \left(0,1\right)$ the increment $J\left(x^{\left(+\right)} \, \left(\cdot ;\tilde{\vartheta },\varepsilon \right)\right)-J\left(\bar{x}\left(\cdot \right)\right)=:\Delta _{\varepsilon }^{\left(+\right)} J\left(\bar{x}\left(\cdot \right);\tilde{\vartheta }\right)$$\left(J\left(x^{\left(-\right)} \, \left(\cdot ;\tilde{\vartheta },\varepsilon \right)\right)-J\left(\bar{x}\left(\cdot \right)\right)=:\Delta _{\varepsilon }^{\left(-\right)} J\left(\bar{x}\left(\cdot \right);\tilde{\vartheta }\right)\right)$ corresponding to the variation \eqref{GrindEQ__2_34_} ((2.34)), is represented in the form
\[\Delta _{\varepsilon }^{\left(+\right)} J\, \left(\bar{x}\, \left(\cdot \right);\, \tilde{\vartheta }\right)=\varepsilon ^{2} E\left(\bar{L}\right)\, \left(\theta +,\xi \right)+\frac{\varepsilon ^{3} }{2\left(1-\varepsilon \right)} \xi ^{T} \bar{L}_{\dot{x}\dot{x}} \left(\theta +\right)\xi +\] 
\begin{equation} \label{GrindEQ__2_36_} 
+\frac{1}{2} \varepsilon ^{4} \left[K\left(\bar{L}\right)\left(\theta +,\varepsilon ,\, \xi \right)-\frac{1}{3\left(\varepsilon -1\right)^{2} } \left(\xi ^{T} \bar{L}_{\dot{x}\dot{x}} \left(\theta +\right)\xi \right)_{\dot{x}}^{T} \xi \right]+o\left(\varepsilon ^{4} \right) 
\end{equation} 
\[\left(\Delta _{\varepsilon }^{\left(-\right)} J\, \left(\bar{x}\, \left(\cdot \right);\, \tilde{\vartheta }\right)=\varepsilon ^{2} E\left(\bar{L}\right)\, \left(\theta -,\xi \right)+\frac{\varepsilon ^{3} }{2\left(1-\varepsilon \right)} \xi ^{T} \bar{L}_{\dot{x}\dot{x}} \left(\theta -\right)\xi -\right. \] 
\begin{equation} \label{GrindEQ__2_37_} 
\left. -\frac{1}{2} \varepsilon ^{4} \left[K\left(\bar{L}\right)\left(\theta -,\varepsilon ,\, \xi \right)+\frac{1}{3\left(\varepsilon -1\right)^{2} } \left(\xi ^{T} \bar{L}_{\dot{x}\dot{x}} \left(\theta -\right)\xi \right)_{\dot{x}}^{T} \xi \right]+o\left(\varepsilon ^{4} \right)\right).   
\end{equation} 
Here $E(\bar{L})(\theta ,\xi)$ is determined by \eqref{GrindEQ__1_6_} and the function $K(\bar{L})(\theta ,\varepsilon ,\xi)$  has the form 
\[K\left(\bar{L}\right)\left(\theta ,\, \varepsilon ,\xi \right)=\xi ^{T} \left[\bar{L}_{x} \left(\theta ,\xi \right)-\bar{L}_{x} \left(\theta \right)-\bar{L}_{x\dot{x}} \left(\theta \right)\xi \right]+\frac{d}{dt} E\left(\bar{L}\right)\left(\theta ,\xi \right)+\] 
\begin{equation} \label{GrindEQ__2_38_} 
+\frac{1+\varepsilon }{2\left(1-\varepsilon \right)} \frac{d}{dt} \xi ^{T} \bar{L}_{\dot{x}\dot{x}} \left(\theta \right)\xi ,      
\end{equation} 
where $\left(\xi ^{T} \bar{L}_{\dot{x}\dot{x}} \left(\theta \right)\xi \right)_{\dot{x}} :=\left. \left(\xi ^{T} L_{\dot{x}\dot{x}} \left(t,\, x,\, \dot{x}\right)\xi \right)_{\dot{x}} \right|_{x=\bar{x}\left(t\right)} $is a derivative of the function $\xi ^{T} L_{\dot{x}\dot{x}} \left(t,\, x,\, \dot{x}\right)\xi $ with respect to the variable $\dot{x}$, calculated along the extremal $\bar{x}(\cdot)$.

\textbf{Proof. }At first, we prove the validity of expansion formula \eqref{GrindEQ__2_36_}, i.e. we find expansion with respect to $\varepsilon $ of the increment $J\left(x^{\left(+\right)} \, \left(\cdot ;\tilde{\vartheta },\varepsilon \right)\right)-J\left(\bar{x}\left(\cdot \right)\right)=:\Delta _{\varepsilon }^{\left(+\right)} J\left(\bar{x}\left(\cdot \right),\tilde{\vartheta }\right)$, where $\bar{x}(\cdot)$ is an extremal in problem \eqref{GrindEQ__1_1_}, \eqref{GrindEQ__1_2_}, the admissible function $x^{\left(+\right)} (\cdot ;\tilde{\vartheta },\varepsilon)$ is determined by \eqref{GrindEQ__2_34_}.  

Substituting $\lambda =\varepsilon \in \left(0,\, \bar{\varepsilon }\right]\bigcap \left(0,1\right)$, i.e. $\vartheta =\tilde{\vartheta }$ in   \eqref{GrindEQ__2_17_}-\eqref{GrindEQ__2_19_}, we calculate the increment $\Delta _{\varepsilon }^{\left(+\right)} J\left(\bar{x}\left(\cdot \right),\tilde{\vartheta }\right)$ with accuracy $o(\varepsilon ^{4} )$. Then, allowing for \eqref{GrindEQ__1_6_} and \eqref{GrindEQ__2_34_}, we have 
\begin{equation} \label{GrindEQ__2_39_} 
\Delta _{\varepsilon }^{\left(+\right)} J\left(\bar{x}\left(\cdot \right);\tilde{\vartheta }\right)=J_{1}^{\left(+\right)} \left(\varepsilon ,\, \tilde{\vartheta }\right)+J_{2}^{\left(+\right)} \left(\varepsilon ,\, \tilde{\vartheta }\right),\, \, \varepsilon \in \left(0,\, \bar{\varepsilon }\right]\bigcap \left(0,1\right),   
\end{equation} 
where
\begin{equation} \label{GrindEQ__2_40_} 
J_{1}^{\left(+\right)} \left(\varepsilon ,\, \tilde{\vartheta }\right)=\int _{\theta }^{\theta +\varepsilon ^{2} }E\left(\bar{L}\right)\, \left(t,\, \xi \right)dt+ \int _{\theta +\varepsilon ^{2} }^{\theta +\varepsilon }E\left(\bar{L}\right)\left(t,\frac{\varepsilon }{\varepsilon -1} \, \xi \right)dt ,     
\end{equation} 
\[J_{2}^{\left(+\right)} \left(\varepsilon ,\, \tilde{\vartheta }\right)=\int _{\theta }^{\theta +\varepsilon ^{2} }\, \, \, \left[L\left(t,\bar{x}\left(t\right)+\left(t-\theta \right)\xi ,\, \dot{\bar{x}}\left(t\right)+\xi \right)-L\left(t,\bar{x}\left(t\right),\, \dot{\bar{x}}\left(t\right)+\, \xi \right)-\left(t-\theta \right)\, \bar{L}_{x}^{T} \left(t\right)\, \xi \right] \, dt+\] 
\[+\int _{\theta +\varepsilon ^{2} }^{\theta +\varepsilon }\left[L\left(t,\, \bar{x}\left(t\right)+\frac{\varepsilon }{\varepsilon -1} \left(t-\theta -\varepsilon \right)\xi ,\, \dot{\bar{x}}\left(t\right)+\frac{\varepsilon }{\varepsilon -1} \xi \right)\right. \left. -L\left(t,\bar{x}\left(t\right),\dot{\bar{x}}\left(t\right)+\frac{\varepsilon }{\varepsilon -1} \xi \right)\right] dt-\] 
\begin{equation} \label{GrindEQ__2_41_} 
-\frac{\varepsilon }{\varepsilon -1} \int _{\theta +\varepsilon ^{2} }^{\theta +\varepsilon }\left(t-\theta -\varepsilon \right)\bar{L}_{x}^{T} \left(t\right)\xi dt .     
\end{equation} 
Considering the assumption on the smoothness  of functions $L\left(\cdot \right)$, $L_{\dot{x}} \left(\cdot \right)$ and $\bar{x}\left(\cdot \right)$, we apply the Taylor formula. Then from \eqref{GrindEQ__2_40_} and \eqref{GrindEQ__2_41_}, allowing for notations \eqref{GrindEQ__1_6_} and \eqref{GrindEQ__2_9_}, we get
\[J_{1}^{\left(+\right)} \left(\varepsilon ,\, \tilde{\vartheta }\right)=\int _{\theta }^{\theta +\varepsilon ^{2} }\, \left[E\left(\bar{L}\right)\, \left(\theta +,\, \xi \right)+\left(t-\theta \right)\frac{d}{dt} \, E\left(\bar{L}\right)\, \left(\theta +,\, \xi \right)+o\, \left(t-\theta \right)\right]dt+ \] 
\[+\int _{\theta +\varepsilon ^{2} }^{\theta +\varepsilon }\left[\frac{\varepsilon ^{2} }{2\left(\varepsilon -1\right)^{2} } \xi ^{T} \bar{L}_{\dot{x}\dot{x}} \left(t\right)\, \xi +\frac{1}{6} \frac{\varepsilon ^{3} }{\left(\varepsilon -1\right)^{3} } \left(\xi ^{T} L_{\dot{x}\dot{x}} \left(t\right)\xi \right)\, _{\dot{x}}^{T} \, \xi +o\left(\varepsilon ^{3} ;t\right)\right]dt= \] 
\[=\varepsilon ^{2} E\left(\bar{L}\right)\, \left(\theta +,\, \xi \right)+\frac{\varepsilon ^{3} }{2\left(1-\varepsilon \right)} \xi ^{T} \bar{L}_{\dot{x}\dot{x}} \left(\theta +\right)\, \xi +\frac{\varepsilon ^{4} }{2} \left[\frac{d}{dt} E\left(\bar{L}\right)\, \left(\theta +,\, \xi \right)+\right. \] 
\begin{equation} \label{GrindEQ__2_42_} 
+\frac{1+\varepsilon }{2\left(1-\varepsilon \right)} \frac{d}{dt} \xi ^{T} \bar{L}_{\dot{x}\dot{x}} \left(\theta +\right)\, \xi -\left. \frac{1}{3\left(1-\varepsilon \right)^{2} } \left(\xi ^{T} \bar{L}_{\dot{x}\dot{x}} \left(\theta +\right)\, \xi \right)\, _{\dot{x}}^{T} \xi \right]+o\left(\varepsilon ^{4} \right),   
\end{equation} 
\[J_{2}^{\left(+\right)} \left(\varepsilon ,\, \tilde{\vartheta }\right)\, =\int _{\theta }^{\theta +\varepsilon ^{2} }\, \, \, \left[\left(t-\theta \right)\, \left(\bar{L}_{x}^{T} \left(t,\xi \right)-\bar{L}_{x}^{T} \left(t\right)\right)\xi +o\, \left(t-\theta \right)\, \right] \, \, dt+\] 
\[+\int _{\theta +\varepsilon ^{2} }^{\theta +\varepsilon }\, \frac{\varepsilon }{\varepsilon -1} \left(t-\theta -\varepsilon \right)\, \, \left[\left. \left(\bar{L}_{x}^{T} \left(t,\, \frac{\varepsilon }{\varepsilon -1} \xi \right)-\bar{L}_{x}^{T} \left(t\right)\right)\, \xi \right]\right.  dt+\] 
\[+\int _{\theta +\varepsilon ^{2} }^{\theta +\varepsilon }\left[\frac{1}{2} \left(t-\theta -\varepsilon \right)^{2} \frac{\varepsilon ^{2} }{\left(\varepsilon -1\right)^{2} } \xi ^{T} \bar{L}_{xx} \left(t,\, \frac{\varepsilon }{\varepsilon -1} \xi \right)\xi +o\left(\left(\, \frac{\varepsilon \, \left(t-\theta -\varepsilon \right)}{\varepsilon -1} \right)^{2} \right)\right] \, dt=\] 
\begin{equation} \label{GrindEQ__2_43_} 
=\frac{\varepsilon ^{4} }{2} \left[\xi ^{T} \left(\bar{L}_{x} \left(\theta +,\xi \right)-\bar{L}_{x} \left(\theta +\right)\right)-\xi ^{T} \bar{L}_{x\dot{x}} \left(\theta +\right)\, \xi \right]+o\left(\varepsilon ^{4} \right).  
\end{equation} 

Substituting \eqref{GrindEQ__2_42_} and \eqref{GrindEQ__2_43_} in \eqref{GrindEQ__2_39_} and taking into account \eqref{GrindEQ__2_38_}, and also choosing $\varepsilon ^{*} =\bar{\varepsilon }$, we get the validity of the expansion formula \eqref{GrindEQ__2_36_}. 

\noindent The proof of the validity of the increment formula \eqref{GrindEQ__2_37_} is carried out similar to  \eqref{GrindEQ__2_36_} using \eqref{GrindEQ__2_5_}, \eqref{GrindEQ__2_6_} and \eqref{GrindEQ__2_25_}-\eqref{GrindEQ__2_27_}. Namely, substituting $\lambda =\varepsilon \in \left(0,\tilde{\varepsilon }\right]\bigcap \left(0,1\right)$, i.e. $\vartheta =\tilde{\vartheta }$ in \eqref{GrindEQ__2_25_}-\eqref{GrindEQ__2_27_}, we calculate the increment $J\left(x^{\left(-\right)} \, \left(\cdot ;\tilde{\vartheta },\varepsilon \right)\right)-J\left(\bar{x}\left(\cdot \right)\right)=:\Delta _{\varepsilon }^{\left(-\right)} J\left(\bar{x}\left(\cdot \right);\tilde{\vartheta }\right)$ with accuracy $o\left(\varepsilon ^{4} \right)$, where the number $\tilde{\varepsilon }>0$ is determined above (see (2.5)), while the admissible function $x^{\left(-\right)} (\cdot ;\tilde{\vartheta },\varepsilon)$ by \eqref{GrindEQ__2_35_}. 

Then we have
\begin{equation} \label{GrindEQ__2_44_} 
\Delta _{\varepsilon }^{\left(-\right)} J\left(\, \bar{x}\, \left(\cdot \right),\, \tilde{\vartheta }\right)=J_{1}^{\left(-\right)} \left(\varepsilon ,\, \tilde{\vartheta }\right)+J_{2}^{\left(-\right)}(\varepsilon ,\, \tilde{\vartheta }), \,\, \varepsilon \in \left(0,\tilde{\varepsilon }\right]\bigcap \left(0,1\right),   
\end{equation} 
where
\[J_{1}^{\left(-\right)} \left(\varepsilon ,\, \tilde{\vartheta }\right)=\int _{\theta -\varepsilon ^{2} }^{\theta }\, E\left(\bar{L}\right)\, \left(t,\, \xi \right)dt+ \int _{\theta -\varepsilon }^{\theta -\varepsilon ^{2} }E\left(\bar{L}\right)\, \left(t,\, \frac{\varepsilon }{\varepsilon -1} \xi \right)dt ,   \] 
\[J_{2}^{\left(-\right)} \left(\varepsilon ,\, \tilde{\vartheta }\right)=\int _{\theta -\varepsilon ^{2} }^{\theta }\left[L\, \left(t,\, \bar{x}\left(t\right)+\left(t-\theta \right)\, \xi ,\dot{\bar{x}}\left(t\right)+\xi \, \right)-\, \right.  \left. L\, \left(t,\, \bar{x}\left(t\right),\, \dot{\bar{x}}\left(t\right)+\xi -\left(t-\theta \right)\, \bar{L}_{x}^{T} \left(t\right)\, \xi \right)\, \right]\, dt+\] 
\[+\int _{\theta -\varepsilon }^{\, \theta -\varepsilon ^{2} }\left[L\left(t,\, \bar{x}\left(t\right)+\frac{\varepsilon }{\varepsilon -1} \right. \, \left(t-\theta +\varepsilon \right)\right.  \left. \, \xi ,\dot{\bar{x}}\left(t\right)\left. +\frac{\varepsilon }{\varepsilon -1} \xi \right)-L\left(t,\, \bar{x}\left(t\right),\, \dot{\bar{x}}\left(t\right)+\frac{\varepsilon }{\varepsilon -1} \xi \right)\right]dt-\] 
\[-\frac{\varepsilon }{\varepsilon -1} \int _{\theta -\varepsilon }^{\theta -\varepsilon ^{2} }\left(t-\theta +\varepsilon \right)\, \bar{L}_{x}^{T} \left(t\right)\xi  \, dt.    \] 
Applying the Taylor formula, allowing for notations \eqref{GrindEQ__1_6_}, \eqref{GrindEQ__2_9_} and \eqref{GrindEQ__2_10_}, for the integrals $J_{1}^{\left(-\right)} \left(\varepsilon ,\, \tilde{\vartheta }\right)$ and $J_{2}^{\left(-\right)} \left(\varepsilon ,\, \tilde{\vartheta }\right)$ we have the following expansions:
\[J_{1}^{\left(-\right)} \left(\varepsilon ,\, \tilde{\vartheta }\right)=\int _{\theta -\varepsilon ^{2} }^{\theta }\left[E\, \left(\bar{L}\right)\, \left(\theta -,\, \xi \right)+\left(t-\theta \right)\, \frac{d}{dt} E\, \left(\bar{L}\right)\, \left(\theta -,\, \xi \right)+o\, \left(t-\theta \right)\right] \, \, dt+\] 
\[+\int _{\theta -\varepsilon }^{\theta -\varepsilon ^{2} }\left[\frac{\varepsilon ^{2} }{2\left(\varepsilon -1\right)^{2} } \xi ^{T} \bar{L}_{\dot{x}\dot{x}} \left(t\right)\xi +\frac{1}{6} \frac{\varepsilon ^{3} }{\left(\varepsilon -1\right)^{3} } \, \left(\xi ^{T} \bar{L}_{\dot{x}\dot{x}} \left(t\right)\, \xi \right)_{\dot{x}}^{T} \xi +o\, \left(\varepsilon ^{3} ;\, t\right)\right] \, dt=\] 
\[=\varepsilon ^{2} E\, \left(\bar{L}\right)\, \left(\theta -,\, \xi \right)+\frac{\varepsilon ^{3} }{2\left(1-\varepsilon \right)} \xi ^{T} \bar{L}_{\dot{x}\dot{x}} \, \left(\theta -\right)\, \xi -\frac{1}{2} \, \varepsilon ^{4} \left[\frac{d}{dt} E\, \left(\bar{L}\right)\, \left(\theta -,\, \xi \right)+\right. \] 
\begin{equation} \label{GrindEQ__2_45_} 
+\frac{1+\varepsilon }{2\left(1-\varepsilon \right)} \frac{d}{dt} \left. \xi ^{T} \bar{L}_{\dot{x}\dot{x}} \left(\theta -\right)\xi +\frac{1}{3\left(1-\varepsilon \right)^{2} } \, \left(\xi ^{T} \bar{L}_{\dot{x}\dot{x}} \left(\theta -\right)\, \xi \right)_{\dot{x}}^{\, T} \xi \right]+o\left(\varepsilon ^{4} \right),   
\end{equation} 
\[J_{2}^{\left(-\right)} \left(\varepsilon ,\, \tilde{\vartheta }\right)=\int _{\theta -\varepsilon ^{2} }^{\theta }\left[\, \left(t-\theta \right)\, \Delta \bar{L}_{x}^{T} \left(t,\, \xi \right)\, \xi +o\, \left(t-\theta \right)\, \right] \, dt+\] 
\[+\int _{\theta -\varepsilon }^{\theta -\varepsilon ^{2} }\left[\frac{\varepsilon \left(t-\theta +\varepsilon \right)}{\varepsilon -1} \Delta \bar{L}_{x}^{T} \left(t,\frac{\varepsilon }{\varepsilon -1} \, \xi \right)\xi +\frac{\varepsilon ^{2} \left(t-\theta +\varepsilon \right)^{2} }{2\left(\varepsilon -1\right)^{2} } \xi ^{T} \bar{L}_{xx} \left(t,\, \frac{\varepsilon }{\varepsilon -1} \xi \right)\xi \right]dt+ \] 
\[+\int _{\theta -\varepsilon }^{\theta -\varepsilon ^{2} }o\left(\left(\frac{\varepsilon \left(t-\theta +\varepsilon \right)}{\varepsilon -1} \right)^{2} \right)dt =-\frac{\varepsilon ^{4} }{2} \left[\xi ^{T} \left(\bar{L}_{x} \left(\theta -,\, \xi \right)-\bar{L}_{x} \left(\theta -\right)\right)-\right. \] 
\begin{equation} \label{GrindEQ__2_46_} 
\left. -\xi ^{T} \bar{L}_{x\dot{x}} \left(\theta -\right)\, \xi \, \, \, \right]+o\left(\varepsilon ^{4} \right).      
\end{equation} 

Consequently, substituting  \eqref{GrindEQ__2_45_} and \eqref{GrindEQ__2_46_} in \eqref{GrindEQ__2_44_}, and taking into account \eqref{GrindEQ__2_38_} and choosing $\varepsilon ^{*} =\tilde{\varepsilon }$, we get the validity of the increment formula \eqref{GrindEQ__2_37_}. Thereby, Proposition 2.3 is proved. 

\textbf{Remark 2.1. }Based on Proposition 2.3, we confirm that various methods for choosing the parameter $\lambda $, as the function $\varepsilon $, allows to get a new increment formula of functional in problem \eqref{GrindEQ__1_1_}, \eqref{GrindEQ__1_2_}.

\section{Necessary conditions for a minimum in the presence of various degenerations at a point}

In this section, using the results of the previous section we get various necessary conditions for a strong and weak local minimum with the degeneration of the Weierstrass condition and also with degenerations of Weierstrass and Legendre conditions simultaneously.   

\textbf{Theorem 3.1. }Let the functions $L\left(\cdot \right)$ and $L_{\dot{x}} \left(\cdot \right)$ be continuously differentiable in totality of variables, and the admissible function $\bar{x}\left(\cdot \right)$ be a strong local minimum in problem \eqref{GrindEQ__1_1_}, \eqref{GrindEQ__1_2_}. Then:

(i) if at the point $\theta \in \left[t_{0} ,\, t_{1} \right)$ $\left(\theta \in \left(t_{0} ,\, t_{1} \right]\, \right)$ the function $\bar{x}\left(\cdot \right)$ is twice differentiable on the right (left) in semi-neighborhood $\left[\theta ,\, \theta +\alpha \right)\subset I$  $\left(\left(\theta -\alpha ,\, \theta \right]\subset I\right)$ of the point $\theta $, and also along it for the vectors $\eta \ne 0$ and $\left(\bar{\lambda }-1\right)^{-1} \bar{\lambda }\eta $, where $\bar{\lambda }\in \left(0,\, 1\right)$, the Weierstrass condition degenerates at the point $\theta $ on the right (left) i.e. the following equality holds:
\begin{equation} \label{GrindEQ__3_1_} 
E(\bar{L})(\theta +,\eta)=E\left(\bar{L}\right)(\theta +,\left(\bar{\lambda }-1\right)^{-1} \bar{\lambda }\eta)=0 
\end{equation} 
\begin{equation} \label{GrindEQ__3_2_} 
(E(\bar{L})(\theta -,\eta)=E(\bar{L})(\theta -,(\bar{\lambda}-1)^{-1} \bar{\lambda }\eta)=0),   
\end{equation} 
then the following inequality is fulfilled:
\begin{equation} \label{GrindEQ__3_3_} 
\bar{\lambda }M_{1} \left(\bar{L}_{x} \right)(\theta +,\bar{\lambda },\eta)+\frac{d}{dt} Q_{2} \left(\bar{L}\right)(\theta +,\bar{\lambda },\eta)\ge 0 
\end{equation} 
\begin{equation} \label{GrindEQ__3_4_} 
\left(\bar{\lambda }M_{1} \left(\bar{L}_{x} \right)(\theta -,\bar{\lambda },\eta)+\frac{d}{dt} Q_{2} \left(\bar{L}\right)(\theta -,\bar{\lambda },\eta)\le 0\right),       
\end{equation} 
where $E(\bar{L})(\cdot)$, $Q_{2} (\bar{L})(\cdot)$ and $M_{1} (\bar{L}_{x} )(\cdot)$ are determined by \eqref{GrindEQ__1_6_}, \eqref{GrindEQ__2_11_}, \eqref{GrindEQ__2_12_} allowing for \eqref{GrindEQ__2_9_} and \eqref{GrindEQ__2_10_};  

(ii) if at the point $\theta \in \left(t_{0} ,\, t_{1} \right)$ the function $\bar{x}\left(\cdot \right)$ is twice differentiable, furthermore, along it for the vectors $\eta \ne 0$ and $\left(\lambda -1\right)^{-1} \bar{\lambda }\eta $, where $\bar{\lambda }\in \left(0,1\right)$, the Weierstrass condition degenerates at the point $\theta $, i.e. we have the following equalities 
\begin{equation} \label{GrindEQ__3_5_} 
E\left(\bar{L}\right)(\theta ,\, \eta)=E\left(\bar{L}\right)\left(\theta ,\, \left(\bar{\lambda }-1\right)^{-1} \bar{\lambda }\eta \right)=0,      
\end{equation} 
then the following equality is fulfilled:
\begin{equation} \label{GrindEQ__3_6_} 
\bar{\lambda }\Delta \bar{L}_{x}^{T} \left(\theta ,\, \eta \right)\eta +\left(1-\bar{\lambda }\right)\Delta \bar{L}_{x}^{T} \left(\theta ,\, \left(\bar{\lambda }-1\right)^{-1} \bar{\lambda }\eta \right)\, \eta =0,    
\end{equation} 
where $\Delta \bar{L}_{x} \left(\cdot \right)$ is determined from \eqref{GrindEQ__2_10_} allowing for \eqref{GrindEQ__2_9_}.  

 \textbf{Proof. }At first we prove part (i) of Theorem 3.1, i.e. the validity of inequality \eqref{GrindEQ__3_3_}  (\eqref{GrindEQ__3_4_}).  We use Proposition 2.2 (this is possible by the conditions of Theorem 3.1). Assume $\xi =\eta $ and $\lambda =\bar{\lambda }$, i.e. $\vartheta =\bar{\vartheta }=\left(\theta ,\, \bar{\lambda },\, \eta \right)$ in statement \eqref{GrindEQ__2_32_} ((2.33)) of Proposition  2.2. Then by virtue of assumption \eqref{GrindEQ__3_1_} ((3.2)) and denotation \eqref{GrindEQ__2_11_}, allowing for \eqref{GrindEQ__2_15_}, the increment formula \eqref{GrindEQ__2_32_} (\eqref{GrindEQ__2_33_}) takes the form
\begin{equation} \label{GrindEQ__3_7_} 
\Delta _{\varepsilon }^{\left(+\right)} J\left(\bar{x}\left(\cdot \right),\, \bar{\vartheta }\right)=\frac{1}{2} \varepsilon ^{2} \left[\bar{\lambda }M_{1} \left(\bar{L}_{x} \right)\, \left(\theta +,\, \bar{\lambda },\eta \right)+\frac{d}{dt} Q_{2} \left(\bar{L}\right)\, \left(\theta +,\, \lambda ,\eta \right)\right]
+o\left(\varepsilon ^{2} \right),\, \varepsilon \in \left(0,\, \varepsilon ^{*} \right] 
\end{equation} 
\begin{equation} \label{GrindEQ__3_8_} 
\left(\Delta _{\varepsilon }^{\left(-\right)} J\left(\bar{x}\left(\cdot \right),\, \bar{\vartheta }\right)=-\frac{1}{2} \varepsilon ^{2} \left[\bar{\lambda }M_{1} \left(\bar{L}_{x} \right)\, \left(\theta -,\, \bar{\lambda },\eta \right)+\frac{d}{dt} Q_{2} \left(\bar{L}\right)\, \left(\theta -,\, \lambda ,\eta \right)\right]+\right.    \left. o\left(\varepsilon ^{2} \right),\, \varepsilon \in \left(0,\, \varepsilon ^{*} \right]\, \right).      
\end{equation} 

Since the admissible function $\bar{x}\left(\cdot \right)$ is a strong local minimum in problem \eqref{GrindEQ__1_1_}, \eqref{GrindEQ__1_2_}, then 

\[\varepsilon ^{-2} \Delta _{\varepsilon }^{\left(+\right)} \, J\left(\bar{x}\left(\cdot \right),\, \bar{\vartheta }\right)\ge 0 \,\,\, \left(\varepsilon ^{-2} \Delta _{\varepsilon }^{\left(-\right)} \, J\left(\bar{x}\left(\cdot \right),\, \bar{\vartheta }\right)\ge 0\right)\]  for all $\varepsilon \in \left(0,\, \varepsilon ^{*} \right]$. Therefore, allowing for \eqref{GrindEQ__3_7_} ((3.8)), in the last inequality we pass to the limit as $\varepsilon \to +0$. Then we get the validity of the sought-for inequality \eqref{GrindEQ__3_3_} ((3.4)), i.e. part (i) of Theorem 3.1 is proved.

We now prove part (ii)  of Theorem 3.1, i.e. the validity of the equality \eqref{GrindEQ__3_6_}. Since $\theta \in \left(t_{0} ,\, t_{1} \right)$, and in addition the function $\bar{x}\left(\cdot \right)$ is a strong local minimum in problem \eqref{GrindEQ__1_1_}, \eqref{GrindEQ__1_2_}, and is twice differentiable at the point $\theta $, then from inequality \eqref{GrindEQ__3_3_} and \eqref{GrindEQ__3_4_}, taking into account \eqref{GrindEQ__2_11_} and \eqref{GrindEQ__2_12_}, we get the validity of the following equality 
\[\bar{\lambda }\, \left[\, \, \bar{\lambda }\Delta \bar{L}_{x}^{T} \left(\theta ,\, \eta \right)\eta +\left(1-\bar{\lambda }\right)\Delta \bar{L}_{x}^{T} \left(\theta ,\, \left(\bar{\lambda }-1\right)^{-1} \bar{\lambda }\eta \right)\, \eta \, \right]\, +\] 
\begin{equation} \label{GrindEQ__3_9_} 
+\frac{d}{dt} \left[\bar{\lambda }^{2} E\left(\bar{L}\right)\, \left(\theta ,\, \eta \right)\eta +\left(1-\bar{\lambda }^{2} \right)E\left(\bar{L}\right)\, \left(\theta ,\, \left(\bar{\lambda }-1\right)^{-1} \bar{\lambda }\eta \right)\, \right]=0.   
\end{equation} 

Since $\bar{x}\left(\cdot \right)$ is a strong local minimum in problem \eqref{GrindEQ__1_1_}, \eqref{GrindEQ__1_2_}, by virtue of the Weierstrass condition \eqref{GrindEQ__1_5_}, allowing for assumption \eqref{GrindEQ__3_5_}, the functions $E(\bar{L})(t,\, \eta),\, t\in I$ and $E\left(\bar{L}\right)\, \left(t,\, \left(\bar{\lambda }-1\right)^{-1} \bar{\lambda }\eta \right),t\in I$, with respect to the variable $t$ attain a minimum at the point $\theta \in \left(t_{0} ,\, t_{1} \right)$. Then, taking into account the smoothness of the functions $L\left(\cdot \right)$, $L_{\dot{x}} \left(\cdot \right)$ and $\bar{x}\left(\cdot \right)$, by the Fermat theorem [11, p.15] we have $\frac{d}{dt} E\left(\bar{L}\right)\, \left(\theta ,\, \eta \right)=\frac{d}{dt} E\left(\bar{L}\right)\, \left(\theta ,\, \left(\bar{\lambda }-1\right)^{-1} \bar{\lambda }\eta \right)=0$. Consequently, from \eqref{GrindEQ__3_9_}, allowing for the last equalities, we get the proof of the equality \eqref{GrindEQ__3_6_}, i.e. part (ii) of Theorem 3.1 is proved.  Thus, Theorem 3.1 is completely proved. 

Now we consider the case when condition for a minimum \eqref{GrindEQ__3_3_} ((3.4)) also degenerates i.e. we have the equalities 
\begin{equation} \label{GrindEQ__3_10_} 
W(\bar{L})(\theta +,\bar{\lambda },\, \eta)=0 
\end{equation} 

\begin{equation}\label{GrindEQ__3_11_} 
(W(\bar{L})(\theta-, \bar{\lambda}, \eta)=0),
\end{equation} 

where $W(\bar{L})(\cdot ,\bar{\lambda },\, \eta)$ is determined from \eqref{GrindEQ__2_15_}. In this case the following theorem is valid.

 \textbf{Theorem 3.2. }Let the functions $L\left(\cdot \right)$ and $L_{\dot{x}} \left(\cdot \right)$ be twice continuously differentiable in totality of variables, and the admissible function $\bar{x}\left(\cdot \right)$ be a strong local minimum in problem \eqref{GrindEQ__1_1_}, \eqref{GrindEQ__1_2_}. Then:

(i)\textbf{ }if at the point $\theta \in \left[t_{0} ,\, t_{1} \right)$ $\left(\theta \in \left(t_{0} ,\, t_{1} \right]\right)$ the function $\bar{x}\left(\cdot \right)$ is triply differentiable on the right (left) in semi-neighborhood $\left[\theta ,\, \theta +\alpha \right)\subset I$ $\left(\, \left(\theta -\alpha ,\, \theta \right]\subset I\right)$ of the point $\theta $, and also along it for the vectors $\eta \ne 0$ and  $\left(\bar{\lambda }-1\right)^{-1} \bar{\lambda} \eta $, where $\bar{\lambda }\in \left(0,1\right)$, conditions \eqref{GrindEQ__3_1_} ((3.2)) and \eqref{GrindEQ__3_10_} ((3.11)) are fulfilled, i.e. the Weierstrass condition and the condition for a minimum \eqref{GrindEQ__3_3_} ((3.4)) degenerate at the point $\theta $ on the right (left), then the following inequality is fulfilled:
\begin{equation} \label{GrindEQ__3_12_} 
G\left(\bar{L}\right)\, \left(\theta +,\, \bar{\lambda },\, \eta \right)\ge 0 
\end{equation} 
\begin{equation} \label{GrindEQ__3_13_} 
\left(G\left(\bar{L}\right)\, \left(\theta -,\, \bar{\lambda },\, \eta \right)\ge 0\right);     
\end{equation} 

 (ii) if at the point $\theta \in \left(t_{0} ,t_{1} \right)$ the admissible function $\bar{x}\left(\cdot \right)$ is triply differentiable, furthermore, along it for the vectors $\eta \ne 0$ and $\left(\bar{\lambda }-1\right)^{-1} \bar{\lambda }\eta $, where $\bar{\lambda }\in \left(0,1\right)$, the Weierstrass condition degenerates at the point $\theta $, i.e. the equality \eqref{GrindEQ__3_5_} holds, then the following inequality is fulfilled
\begin{equation} \label{GrindEQ__3_14_} 
G\left(\bar{L}\right)\left(\theta ,\, \bar{\lambda },\, \eta \right)\ge 0,      
\end{equation} 
where $G\left(\bar{L}\right)\left(\cdot ,\bar{\lambda },\eta \right)$ is determined from \eqref{GrindEQ__2_16_}, allowing for \eqref{GrindEQ__2_9_}-\eqref{GrindEQ__2_12_}.

 \textbf{Proof. }Prove part (i) of Theorem 3.2, i.e. the validity of inequality \eqref{GrindEQ__3_12_} (\eqref{GrindEQ__3_13_}). Use Proposition 2.1 (it is possible by virtue of condition of Theorem 3.2). Assume $\xi =\eta $ and $\lambda =\bar{\lambda }$, i.e. $\vartheta =\bar{\vartheta }=\left(\theta ,\, \bar{\lambda },\eta \right)$ in statement \eqref{GrindEQ__2_13_} (\eqref{GrindEQ__2_14_} ) of Proposition 2.1. Since the function $\bar{x}\left(\cdot \right)$ is a strong local minimum in problem \eqref{GrindEQ__1_1_}, \eqref{GrindEQ__1_2_}, then $\Delta _{\varepsilon }^{\left(+\right)} J\left(\bar{x}\left(\cdot \right);\bar{\vartheta }\right)\ge 0$ $\left(\Delta _{\varepsilon }^{\left(-\right)} J\left(\bar{x}\left(\cdot \right);\bar{\vartheta }\right)\ge 0\right)$, $\forall \varepsilon \in \left(0,\, \varepsilon ^{*} \right]$. Then by virtue of \eqref{GrindEQ__3_1_} (\eqref{GrindEQ__3_2_}) and \eqref{GrindEQ__3_10_} (\eqref{GrindEQ__3_11_}), allowing for notation \eqref{GrindEQ__2_11_}, the first two summands in expansion formula \eqref{GrindEQ__2_13_} (\eqref{GrindEQ__2_14_}) vanish. Therefore,  dividing the obtained inequality for $\Delta _{\varepsilon }^{\left(+\right)} J\left(\bar{x}\left(\cdot \right);\, \bar{\vartheta }\right)$ $\left(\Delta _{\varepsilon }^{\left(-\right)} J\left(\bar{x}\left(\cdot \right);\, \bar{\vartheta }\right)\right)$ by $\varepsilon ^{3}$ and passing to the limit as $\varepsilon \to +0$, we get the sought-for inequality \eqref{GrindEQ__3_12_} (\eqref{GrindEQ__3_13_}), i.e. part (i) of Theorem 3.2 is proved. 

 Prove part (ii) of Theorem 3.2, i.e. the validity of inequality \eqref{GrindEQ__3_14_}. Considering the assumption of part (ii) of Theorem 3.2, it is easy to get the validity of equality \eqref{GrindEQ__3_9_}, that was obtained when proving part (ii) of Theorem 3.1. By virtue of notations \eqref{GrindEQ__2_11_}, \eqref{GrindEQ__2_12_} and \eqref{GrindEQ__2_15_}, equality \eqref{GrindEQ__3_9_} takes a new form $W\left(\bar{L}\right)(\theta ,\, \bar{\lambda },\, \eta )=0$. Consequently, we confirm that all assumptions of part (i) of Theorem 3.2 are fulfilled, and the function $\bar{x}\left(\cdot \right)$ is triple differentiable at the point $\theta $. Therefore, the proof of part  (ii) of Theorem 3.2 follows from the statement of part (i) of Theorem 3.2. Theorem 3.2 is completely proved. 

We prove the following theorem in the presence of new degenerations.

\textbf{Theorem 3.3. }Let the admissible function $\bar{x}\left(\cdot \right)$ be a strong local minimum in problem \eqref{GrindEQ__1_1_}, \eqref{GrindEQ__1_2_}. Then:

(i) if the integrant $L\left(\cdot \right)$ is continuously differentiable in totality of variables, and it is triply continuously differentiable with respect to the variables $\dot{x}$, furthermore, along the function $\bar{x}\left(\cdot \right)$ for the vector $\eta \ne 0$ the Legendre condition degenerates at the point $\theta $ or at the point $\theta $ on the right (left), i.e. the equalities 
\begin{equation} \label{GrindEQ__3_15_} 
\eta ^{T} \bar{L}_{\dot{x}\dot{x}} \left(\theta \right)\eta =0,\, \, \,\,\,\, \eta ^{T} \bar{L}_{\dot{x}\dot{x}} \left(\theta +\right)\, \eta =0\, \, \, \,\,\,\, \left(\eta ^{T} \bar{L}_{\dot{x}\dot{x}} \left(\theta -\right)\eta =0\right),    
\end{equation} 
hold, then the following equalities are fulfilled: 
\begin{equation} \label{GrindEQ__3_16_} 
\left(\eta ^{T} \bar{L}_{\dot{x}\dot{x}} \left(\theta \right)\eta \right)_{\, \dot{x}}^{T} \eta =0,\, \, \, \left(\eta ^{T} \bar{L}_{\dot{x}\dot{x}} \left(\theta +\right)\, \eta \right)_{\, \dot{x}}^{T} \eta \, =0\, \, \left(\left(\eta ^{T} \bar{L}_{\dot{x}\dot{x}} \left(\theta -\right)\eta \right)_{\, \dot{x}}^{T} \eta =0\right),   
\end{equation} 
where $\theta \in \left(t_{0} ,\, t_{1} \right)\backslash \left\{\tau \right\}$ or $\theta \in \left\{t_{0} \right\}\bigcup \left\{\tau \right\}\, \, \left(\theta \in \left\{\tau \right\}\bigcup \left\{t_{1} \right\}\right)$, moreover $\left\{\tau \right\}$ is the set of angular points of the function $\bar{x}\left(\cdot \right)$, the symbol $\left(\eta ^{T} \bar{L}_{\dot{x}\dot{x}} \left(\cdot \right)\eta \right)\, _{\dot{x}} $ is determined above (see (2.36));  

(ii) if the functions $L\left(\cdot \right)$ and $L_{\dot{x}} \left(\cdot \right)$ are twice differentiable in totality of variables, furthermore, at the point $\theta \in \left[t_{0} ,\, t_{1} \right)$ $\left(\theta \in \left(t_{0} ,\, t_{1} \right]\right)$ the function $\bar{x}\left(\cdot \right)$ is twice differentiable in semi-neighborhood $\left[\theta ,\, \theta +\alpha \right)\subset I$ $\left(\left(\theta -\alpha ,\, \theta \right]\subset I\right)$ of the point $\theta $, and along it for the vector $\eta \ne 0$ the Weierstrass and Legendre conditions degenerate at the point $\theta $ on the right (left), i.e. the equalities 
\begin{equation} \label{GrindEQ__3_17_} 
E\left(\bar{L}\right)\left(\theta +,\, \eta \right)=\eta ^{T} \bar{L}_{\dot{x}\dot{x}} \left(\theta +\right)\, \eta =0 
\end{equation} 
\begin{equation} \label{GrindEQ__3_18_} 
\left(E\left(\bar{L}\right)\left(\theta -,\, \eta \right)=\eta ^{T} \bar{L}_{\dot{x}\dot{x}} \left(\theta -\right)\, \eta =0\right),     
\end{equation} 
hold, then the following inequalities are fulfilled:
\[\eta ^{T} \left(\bar{L}_{x} \left(\theta +,\, \eta \right)-\bar{L}_{x} \left(\theta +\right)-\bar{L}_{x\dot{x}} \left(\theta +\right)\eta \right)+\frac{d}{dt} E\left(\bar{L}\right)\left(\theta +,\, \eta \right)+\] 
\begin{equation} \label{GrindEQ__3_19_} 
+\frac{1}{2} \frac{d}{dt} \eta ^{T} \bar{L}_{\dot{x}\dot{x}} \left(\theta +\right)\eta \ge 0  
\end{equation} 
\[\left(\eta ^{T} \left(\bar{L}_{x} \left(\theta -,\, \eta \right)-\bar{L}_{x} \left(\theta -\right)-\bar{L}_{x\dot{x}} \left(\theta -\right)\eta \right)+\frac{d}{dt} E\left(\bar{L}\right)\left(\theta -,\, \eta \right)+\right. \] 
\begin{equation} \label{GrindEQ__3_20_} 
\left. +\frac{1}{2} \frac{d}{dt} \eta ^{T} \bar{L}_{\dot{x}\dot{x}} \left(\theta -\right)\eta \le 0\right),     
\end{equation} 
where $E\left(\bar{L}\right)\left(\cdot ,\, \eta \right)$ and $\bar{L}_{x} \left(\cdot ,\, \eta \right)$ are determined by \eqref{GrindEQ__1_6_} and \eqref{GrindEQ__2_9_}, respectively;

 (iii) if the functions $L\left(\cdot \right)$ and $L_{\dot{x}} \left(\cdot \right)$ are twice continuously differentiable in totality of variables, furthermore, at the point $\theta \in \left(t_{0} ,\, t_{1} \right)$ the function $\bar{x}\left(\cdot \right)$ is twice differentiable and along it for the vector $\eta \ne 0$ the Weierstrass and  Legendre conditions degenerate at the point $\theta $, i.e. the equalities 
\begin{equation} \label{GrindEQ__3_21_} 
E(\bar{L})(\theta ,\, \eta)=\eta ^{T} \bar{L}_{\dot{x}\dot{x}} \left(\theta \right)\, \eta =0    
\end{equation} 
hold, then the following equality is fulfilled:
\begin{equation} \label{GrindEQ__3_22_} 
\, \eta ^{T} \left[L_{x} \left(\theta ,\, \bar{x}\left(\theta \right),\, \dot{\bar{x}}\left(\theta \right)+\eta \right)-\bar{L}_{x} \left(\theta \right)-\bar{L}_{x\dot{x}} \left(\theta \right)\, \eta \right]=0.    
\end{equation} 

\textbf{Proof. }At first we prove statement \eqref{GrindEQ__3_16_}. For that it suffices to show for example the validity of the equality $\left(\eta ^{T} \bar{L}_{\dot{x}\dot{x}} \left(\theta +\right)\eta \right)^T_{\dot{x}}\eta=0$   under the assumption $\eta ^{T} L_{\dot{x}\dot{x}} \left(\theta +\right)\eta =0$, became the validity of other equalities from \eqref{GrindEQ__3_16_} are proved quite similarly. Since $\bar{x}\left(\cdot \right)$ is a strong local minimum in problem \eqref{GrindEQ__1_1_}, \eqref{GrindEQ__1_2_}, then by virtue of \eqref{GrindEQ__1_5_}, \eqref{GrindEQ__1_6_} we have 
\[E\left(\bar{L}\right)\left(\theta +,\, \xi \right)=L\left(\theta ,\, \bar{x}\left(\theta \right),\, \dot{\bar{x}}\left(\theta +\right)+\xi \right)-\bar{L}\left(\theta +\right)-\bar{L}_{\dot{x}} \left(\theta +\right)\xi \ge 0,\, \, \forall \xi \in R^{n} .\] 
Hence, assuming $\xi =\varepsilon \eta $, where $\varepsilon \in \left(-q_{0} ,\, q_{0} \right),\, q_{0} >0$, and considering the smoothness of the integrant $L\left(\cdot \right)$ with respect to the variable $\dot{x}$, by the Taylor formula we get
\[E\left(\bar{L}\right)\left(\theta +,\, \xi \right)=\frac{1}{2} \varepsilon ^{2} \eta ^{T} \bar{L}_{\dot{x}\dot{x}} \left(\theta +\right)\eta +\frac{1}{6} \varepsilon ^{3} \left(\eta ^{T} \bar{L}_{\dot{x}\dot{x}} \left(\theta +\right)\eta \right)_{\dot{x}}^{T} \eta +o\left(\varepsilon ^{3} \right)\ge 0, \] 
\[\forall \varepsilon \in \left(-q_{0} ,\, q_{0} \right) .      \] 

    Since by the assumption $\eta ^{T} L_{\dot{x}\dot{x}} \left(\theta +\right)\eta =0$, then the validity of the equality  $\left(\eta ^{T} L_{\dot{x}\dot{x}} \left(\theta +\right)\eta \right)_{\dot{x}}^{T} \eta =0$ easily follows from the last inequality.

\noindent We now prove the validity of inequality \eqref{GrindEQ__3_19_} (\eqref{GrindEQ__3_20_}), i.e. part (ii) of Theorem 3.3. Let $\bar{x}\left(\cdot \right)$ be a strong local minimum in problem \eqref{GrindEQ__1_1_}, \eqref{GrindEQ__1_2_}. Then by virtue of conditions assumed in part (ii) of Theorem 3.3, we confirm that Proposition 2.3 is valid, namely expansion \eqref{GrindEQ__2_36_} (\eqref{GrindEQ__2_37_}) for the increment $\Delta _{\varepsilon }^{\left(+\right)} J\left(\bar{x}\left(\cdot \right);\tilde{\vartheta }\right)$ $\left(\Delta _{\varepsilon }^{\left(-\right)} J\left(\bar{x}\left(\cdot \right);\tilde{\vartheta }\right)\right)$ and the inequality $\Delta _{\varepsilon }^{\left(+\right)} J\left(\bar{x}\left(\cdot \right);\tilde{\vartheta }\right)\ge 0$ $\left(\Delta _{\varepsilon }^{\left(-\right)} J\left(\bar{x}\left(\cdot \right);\tilde{\vartheta }\right)\ge 0\right)$,  $\forall \varepsilon \in \left(0,\, \varepsilon ^{*} \right]\bigcap \left(0,\, 1\right)$ hold. Assume $\xi =\eta $ in the last inequality and take into account \eqref{GrindEQ__3_17_} ((3.18)), \eqref{GrindEQ__2_36_} (\eqref{GrindEQ__2_37_}), \eqref{GrindEQ__2_38_}, and also statement \eqref{GrindEQ__3_16_} of Theorem 3.3. Then we have:
\[\Delta _{\varepsilon }^{\left(+\right)} J\left(\bar{x}\left(\cdot \right);\tilde{\vartheta }\right)=\frac{1}{2} \varepsilon ^{4} \left[\eta ^{T} \left(\bar{L}_{x} \left(\theta +,\eta \right)-\bar{L}_{x} \left(\theta +\right)-\bar{L}_{x\dot{x}} \left(\theta +\right)\eta \right)+\frac{d}{dt} \right. E\left(\bar{L}\right)\left(\theta +,\, \eta \right)+\] 
\[\left. +\frac{1+\varepsilon }{2\left(1-\varepsilon \right)} \frac{d}{dt} \eta ^{T} \bar{L}_{\dot{x}\dot{x}} \left(\theta +\right)\eta \right]+o\left(\varepsilon ^{4} \right)\ge 0\] 
\[\left(\Delta _{\varepsilon }^{\left(-\right)} J\left(\bar{x}\left(\cdot \right)\tilde{\vartheta }\right)=-\frac{1}{2} \varepsilon ^{4} \left[\eta ^{T} \left(\bar{L}_{x} \left(\theta -,\eta \right)-\bar{L}_{x} \left(\theta -\right)-\bar{L}_{x\dot{x}} \left(\theta -\right)\eta \right)+\frac{d}{dt} \right. E\left(\bar{L}\right)\left(\theta -,\, \eta \right)+\right. \] 
\[\left. \left. +\frac{1+\varepsilon }{2\left(1-\varepsilon \right)} \frac{d}{dt} \eta ^{T} \bar{L}_{\dot{x}\dot{x}} \left(\theta -\right)\eta \right]+o\left(\varepsilon ^{4} \right)\ge 0\right), \forall \varepsilon \in \left(0,\, \varepsilon ^{*} \right]\bigcap \left(0,\, 1\right).  \] 
Dividing the obtained last expression for $\Delta _{\varepsilon }^{\left(+\right)} J\left(\bar{x}\left(\cdot \right);\tilde{\vartheta }\right)$ $\left(\Delta _{\varepsilon }^{\left(-\right)} J\left(\bar{x}\left(\cdot \right);\tilde{\vartheta }\right)\right)$ by $\varepsilon ^{4}$ and passing to the limit as $\varepsilon \to +0$, we get the sought-for inequality \eqref{GrindEQ__3_19_} (\eqref{GrindEQ__3_20_} ). Thus, part (ii) of Theorem 3.3 is proved.

We now prove part (iii) of Theorem 3.3. Since the function $\bar{x}\left(\cdot \right)$ is a strong, local minimum in problem \eqref{GrindEQ__1_1_}, \eqref{GrindEQ__1_2_} and along it \eqref{GrindEQ__3_21_} is fulfilled, we arrive at the conclusion: firstly, assumptions \eqref{GrindEQ__3_17_} and \eqref{GrindEQ__3_18_} are fulfilled, and therefore inequalities \eqref{GrindEQ__3_19_} and \eqref{GrindEQ__3_20_} are valid; secondly, the left hand sides of these inequalities coincide and therefore are equal to zero, i.e. the equality 
\begin{equation} \label{GrindEQ__3_23_} 
\eta ^{T} \left(\bar{L}_{x} \left(\theta ,\, \eta \right)-\bar{L}_{x} \left(\theta \right)-\bar{L}_{x\dot{x}} \left(\theta \right)\eta \right)+\frac{d}{dt} E\left(\bar{L}\right)\left(\theta ,\, \eta \right)+\frac{1}{2} \frac{d}{dt} \eta ^{T} \bar{L}_{\dot{x}\dot{x}} \left(\theta \right)\eta =0 
\end{equation} 
holds; thirdly, by virtue of Weierstrass and Legendre conditions, the functions $E\left(\bar{L}\right)\left(t,\eta \right),$ $t\in I$, and $\eta ^{T} L_{\dot{x}\dot{x}}\left(t,\eta \right)$ $ ,t\in I$, with respect to the variable $t$ obtain a minimum at the point $\theta \in \left(t_{0} ,\, t_{1} \right)$, and therefore, their derivatives with respect to $t$ at the point $\theta $ are equal to zero, i.e.  $\frac{d}{dt} E\left(\bar{L}\right)\left(\theta ,\, \eta \right)=\frac{d}{dt} \eta ^{T} \bar{L}_{\dot{x}\dot{x}} \left(\theta ,\, \eta \right)=0$. Thus, considering the last equality in \eqref{GrindEQ__3_23_}, we get the sought-for equality \eqref{GrindEQ__3_22_}, i.e. part (iii) of Theorem 3.3 is proved. By the same token, Theorem 3.3 is completely proved. 

Continuing the study, below we obtain necessary conditions for a weak local minimum being local modifications of statements of Theorems 3.1-3.3. Namely, we prove the following theorems. \textbf{}

\textbf{Theorem 3.4. }Let the functions $L\left(\cdot \right)$ and $L_{\dot{x}} \left(\cdot \right)$ be continuously differentiable in totality of variables, furthermore, the admissible function $\bar{x}\left(\cdot \right)$ be a weak local minimum in problem \eqref{GrindEQ__1_1_}, \eqref{GrindEQ__1_2_}. Then there exists $\delta >0$ for which the following statements are valid:

(j) if the assumptions of part (i) of Theorem 3.1 are fulfilled, then for every point $\left(\eta ,\, \left(\bar{\lambda }-1\right)^{-1} \bar{\lambda }\eta ,\, \bar{\lambda }\right)\in B_{\delta } \left(0\right)\times B_{\delta } \left(0\right)\times \left(0,\, 1\right)$ satisfying condition \eqref{GrindEQ__3_1_} (\eqref{GrindEQ__3_2_}), the inequality \eqref{GrindEQ__3_3_} (\eqref{GrindEQ__3_4_}) is valid; 

(jj) if the assumptions of part (ii) of Theorem 3.1 are fulfilled, then for every point $\left(\eta ,\, \left(\bar{\lambda }-1\right)^{-1} \bar{\lambda }\eta ,\, \bar{\lambda }\right)\in B_{\delta } \left(0\right)\times B_{\delta } \left(0\right)\times \left(0,\, 1\right)$ satisfying the condition \eqref{GrindEQ__3_5_}, the equality \eqref{GrindEQ__3_6_} is valid. 

\textbf{Proof. }It is clear that by virtue of the assumptions of Theorem 3.4 we have formula \eqref{GrindEQ__3_7_} (\eqref{GrindEQ__3_7_}) obtained for the increment  $\Delta _{\varepsilon }^{\left(+\right)} J\left(\bar{x}\left(\cdot \right),\, \bar{\vartheta }\right)=J\left(x^{\left(+\right)} \left(\cdot ;\, \bar{\vartheta },\, \varepsilon \right)-J\left(\bar{x}\left(\cdot \right)\right)\right)$ $\left(\Delta _{\varepsilon }^{\left(-\right)} J\left(\bar{x}\left(\cdot \right),\, \bar{\vartheta }\right)=J\left(x^{\left(-\right)} \left(\cdot ;\, \bar{\vartheta },\, \varepsilon \right)\right)-J\left(\bar{x}\left(\cdot \right)\right)\right)$, where $x^{\left(+\right)} \left(\cdot ;\, \bar{\vartheta },\, \varepsilon \right)$ $\left(x^{\left(-\right)} \left(\cdot ;\, \bar{\vartheta },\, \varepsilon \right)\right)$ is determined from \eqref{GrindEQ__2_4_}, (\eqref{GrindEQ__2_5_}) allowing for \eqref{GrindEQ__2_1_} (\eqref{GrindEQ__2_6_}) and $\vartheta =\bar{\vartheta }=\left(\theta ,\, \bar{\lambda },\, \eta \right)\in \left[t_{0} ,\, t_{1} \right)\times \left(0,1\right)\times R^{n} \backslash \left\{0\right\}$ $\left(\vartheta =\bar{\vartheta }=\left(\theta ,\, \bar{\lambda },\, \eta \right)\in \left(t_{0} ,\, t_{1} \right]\times \left(0,1\right)\times R^{n} \backslash \left\{0\right\}\right)$. By means of \eqref{GrindEQ__2_1_}, \eqref{GrindEQ__2_2_} and \eqref{GrindEQ__2_4_} ((2.5)-\eqref{GrindEQ__2_7_}), for  $x^{\left(+\right)} \left(\cdot ;\, \bar{\vartheta },\, \varepsilon \right)$ $\left(x^{\left(-\right)} \left(\cdot ;\, \bar{\vartheta },\, \varepsilon \right)\right)$ for all  $\varepsilon \in \left(0,\, \bar{\varepsilon }\right]\bigcap \left(0,\, \tilde{\varepsilon }\right]\bigcap \left(0,\, 1\right]$ the following estimations are valid: 
\[\left\| x^{\left(+\right)} \left(\cdot ;\, \bar{\vartheta },\, \varepsilon \right)-\bar{x}\left(\cdot \right)\right\| _{C\left(I,\, R^{n} \right)} =\left\| h^{\left(+\right)} \left(\cdot ,\, \bar{\vartheta },\, \varepsilon \right)\right\| _{C\left(I,\, R^{n} \right)} \le \left\| \eta \right\| _{R^{n} } ,\] 
\begin{equation} \label{GrindEQ__3_24_} 
\left\| \dot{x}^{\left(+\right)} \left(\cdot ;\, \bar{\vartheta },\, \varepsilon \right)-\dot{\bar{x}}\left(\cdot \right)\right\| _{L_{\infty } \left(I,\, R^{n} \right)} =\left\| \dot{h}^{\left(+\right)} \left(\cdot ;\, \bar{\vartheta },\, \varepsilon \right)\right\| _{L_{\infty } \left(I,\, R^{n} \right)} \le \max \left\{\, \left\| \eta \right\| _{R^{n} } ,\, \left(\bar{\lambda }-1\right)^{-1} \bar{\lambda }\left\| \eta \right\| _{R^{n} } \, \right\} 
\end{equation} 
\[\left(\left\| x^{\left(-\right)} \left(\cdot ;\, \bar{\vartheta },\, \varepsilon \right)-\bar{x}\left(\cdot \right)\right\| _{C\left(I,\, R^{n} \right)} =\left\| h^{\left(-\right)} \left(\cdot ;\bar{\vartheta },\, \varepsilon \right)\right\| _{C\left(I,\, R^{n} \right)} \le \left\| \eta \right\| _{R^{n} } \right. ,  \] 
\begin{equation} \label{GrindEQ__3_25_} 
\left. \left\| \dot{x}^{\left(-\right)} \left(\cdot ;\, \bar{\vartheta },\, \varepsilon \right)-\dot{\bar{x}}\left(\cdot \right)\right\| _{L_{\infty } \left(I,\, R^{n} \right)} =\left\| \dot{h}^{\left(-\right)} \left(\cdot ;\, \bar{\vartheta },\, \varepsilon \right)\right\| _{L_{\infty } \left(I,\, R^{n} \right)} \le \max \left\{\, \left\| \eta \right\| _{R^{n} } ,\, \left(\bar{\lambda }-1\right)^{-1} \bar{\lambda }\left\| \eta \right\| _{R^{n} } \, \right\}\, \right). 
\end{equation} 

\textbf{ }Let the admissible function  $\bar{x}\left(\cdot \right)$ be a weak local minimum in problem \eqref{GrindEQ__1_1_}, \eqref{GrindEQ__1_2_} with $\hat{\delta }$-neighborhood. Then, considering estimations \eqref{GrindEQ__3_24_} (\eqref{GrindEQ__3_25_}), we confirm that for every point $\left(\eta ,\, \left(\bar{\lambda }-1\right)^{-1} \bar{\lambda }\eta ,\, \bar{\lambda }\right)\in B_{\hat{\delta }} \left(0\right)\times B_{\hat{\delta }} \left(0\right)\times \left(0,\, 1\right)$ satisfying condition \eqref{GrindEQ__3_1_} (\eqref{GrindEQ__3_2_}) or \eqref{GrindEQ__3_5_}, the following inequalities are fulfilled:  

\noindent $\Delta _{\varepsilon }^{\left(+\right)} J\left(\bar{x}\left(\cdot \right),\, \bar{\vartheta }\right)\ge 0,\, \forall \varepsilon \in \left(0,\, \bar{\varepsilon }\right]\bigcap \left(0,\, 1\, \right]$ $\left(\Delta _{\varepsilon }^{\left(-\right)} J\left(\bar{x}\left(\cdot \right),\, \bar{\vartheta }\right)\ge 0,\, \forall \varepsilon \in \left(0,\, \bar{\varepsilon }\right]\bigcap \left(0,\, 1\, \right]\right)$ or $\Delta _{\varepsilon }^{\left(+\right)} J\left(\bar{x}\left(\cdot \right),\, \bar{\vartheta }\right)\ge 0$, $\Delta _{\varepsilon }^{\left(-\right)} J\left(\bar{x}\left(\cdot \right),\, \bar{\vartheta }\right)\ge 0$, $\forall \varepsilon \in \left(0,\, \bar{\varepsilon }\right]\bigcap \left(0,\, 1\, \right]$. 

Therefore, allowing for the estimation \eqref{GrindEQ__3_24_}, \eqref{GrindEQ__3_25_}, and also choosing $\delta =\hat{\delta }$, the proof of Theorem 3.4 directly follows from Theorem 3.1. Theorem 3.4 is proved. 

\textbf{ Remark 3.1. }Obviously, if there exists some set $B_{\delta } \left(0\right)\times B_{\delta } \left(0\right)$ that contains no solution of the form $\left(\eta ,\, \left(\bar{\lambda }-1\right)^{-1} \bar{\lambda }\eta \right)\, ,$ where $\eta \ne 0$, $\bar{\lambda }\in \left(0,\, 1\right)$, to each system of equations \eqref{GrindEQ__3_1_}, \eqref{GrindEQ__3_2_} and \eqref{GrindEQ__3_5_}, then Theorem 3.4 is inefficient. 

Using estimations \eqref{GrindEQ__3_24_} and \eqref{GrindEQ__3_25_}, by means of the reason given in the proof of Theorem 3.4, the following statement follows directly from Theorem 3.2. 

\textbf{Theorem 3.5. }Let the functions $L\left(\cdot \right)$ and $L_{\dot{x}} \left(\cdot \right)$ be twice continuously differentiable in totality of variables, furthermore, the admissible function $\bar{x}\left(\cdot \right)$ be a weak local minimum in problem \eqref{GrindEQ__1_1_}, \eqref{GrindEQ__1_2_}. Then there exists $\delta >0$ for which the following statements are valid:\textbf{}

 (j) if the assumptions of part (i) of Theorem 3.2 are fulfilled, then for every point $\left(\eta ,\, \left(\bar{\lambda }-1\right)^{-1} \bar{\lambda }\eta ,\, \bar{\lambda }\right)\in B_{\delta } \left(0\right)\times B_{\delta } \left(0\right)\times \left(0,\, 1\right)$ satisfying conditions  \eqref{GrindEQ__3_1_} (\eqref{GrindEQ__3_2_}) and \eqref{GrindEQ__3_10_} (\eqref{GrindEQ__3_11_}), the inequality \eqref{GrindEQ__3_12_} (\eqref{GrindEQ__3_13_}) is valid;

(jj) if the assumptions of part (ii) of Theorem 3.2 are fulfilled, then for every point $\left(\eta ,\, \left(\bar{\lambda }-1\right)^{-1} \bar{\lambda }\eta ,\, \bar{\lambda }\right)\in B_{\delta } \left(0\right)\times B_{\delta } \left(0\right)\times \left(0,\, 1\right)$ satisfying the condition \eqref{GrindEQ__3_5_}, the inequality \eqref{GrindEQ__3_14_} is valid.

\textbf{Theorem 3.6. }If the assumptions of part (i) of Theorem 3.3 are fulfilled, then the statement \eqref{GrindEQ__3_16_} of Theorem 3.3 is valid also for the admissible function $\bar{x}\left(\cdot \right)$, being a weak local minimum in problem \eqref{GrindEQ__1_1_}, \eqref{GrindEQ__1_2_}. 

\textbf{Proof. }It suffices to show, for example, the validity of the equality  $\left(\eta ^{T} \bar{L}_{\dot{x}\dot{x}} \left(\theta +\right)\eta \right)_{\dot{x}}^{T} \eta =0$ under the assumption
\begin{equation} \label{GrindEQ__3_26_} 
\eta ^{T} \bar{L}_{\dot{x}\dot{x}} \left(\theta +\right)\eta =0,      
\end{equation} 
because other equalities from \eqref{GrindEQ__3_16_} are proved quite similarly. Let the function $\bar{x}\left(\cdot \right)$ be a weak local minimum in problem \eqref{GrindEQ__1_1_}, \eqref{GrindEQ__1_2_} with $\hat{\delta }$-neighborhood and $\eta \in R^{n} $ be an arbitrary fixed vector satisfying the condition \eqref{GrindEQ__3_26_}. We choose a number $\varepsilon _{0} >0$  so that for all $\varepsilon \in \left(-\varepsilon _{0} ,\, \varepsilon _{0} \right)$ the inclusion $\varepsilon \eta \in B_{\hat{\delta }} \left(0\right)$ holds. Assume $\xi =\varepsilon \eta $ in \eqref{GrindEQ__1_8_} and taking into account \eqref{GrindEQ__1_6_}, we apply the Taylor formula. Then by virtue of \eqref{GrindEQ__3_26_} we have 
\[E\left(\bar{L}\right)\left(\theta +,\, \varepsilon \eta \right)=\frac{1}{6} \varepsilon ^{3} \left(\eta ^{T} \bar{L}_{\dot{x}\dot{x}} \left(\theta +\right)\eta \right)_{\dot{x}}^{T} \eta +o\left(\varepsilon ^{3} \right)\ge 0,\, \, \forall \varepsilon \left(-\varepsilon _{0} ,\, \varepsilon _{0} \right].  \] 

Hence we get the validity of the sought-for equality $\left(\eta ^{T} \bar{L}_{\dot{x}\dot{x}} \left(\theta +\right)\eta \right)_{\dot{x}}^{T} \eta =0$. Theorem 3.6 is proved.

\textbf{Theorem 3.7.  }Let the admissible function $\bar{x}\left(t\right)$ be a weak local minimum in problem \eqref{GrindEQ__1_1_}, \eqref{GrindEQ__1_2_}. Then there exists $\delta >0$ for which the following statements are valid: \textbf{}

(j) if the assumptions of part (ii) of Theorem 3.3 are fulfilled, then for every point $\eta \in B_{\delta } \left(0\right)$ satisfying condition \eqref{GrindEQ__3_17_} ((3.18)), the inequality \eqref{GrindEQ__3_19_} ((3.20)) is valid;

(jj) if the assumptions of part (iii) of Theorem 3.3 are fulfilled, then for every point $\eta \in B_{\delta } \left(0\right)$ satisfying condition \eqref{GrindEQ__3_21_}, the inequality \eqref{GrindEQ__3_22_} is valid. 

\textbf{Proof. }We consider the increment formula \eqref{GrindEQ__2_36_} ((2.37)) obtained for  $\Delta _{\varepsilon }^{\left(+\right)} J\left(\bar{x}\left(\cdot \right),\, \tilde{\vartheta }\right)=J\left(x^{\left(+\right)} \left(\cdot ;\tilde{\vartheta },\, \varepsilon \right)\right)-J\left(\bar{x}\left(\cdot \right)\right)$ $\left(\Delta _{\varepsilon }^{\left(-\right)} J\left(\bar{x}\left(\cdot \right),\, \tilde{\vartheta }\right)=J\left(x^{\left(-\right)} \left(\cdot ;\, \tilde{\vartheta },\, \varepsilon \right)\right)-J\left(\bar{x}\left(\cdot \right)\right)\right)$, and used while proving theorem 3.3. 

Here $\left. \tilde{\vartheta }=\vartheta =\left(\theta ,\, \lambda ,\, \xi \right)\right|_{\lambda =\varepsilon } $,$\left. x^{\left(+\right)} \left(\cdot ;\, \tilde{\vartheta },\, \varepsilon \right)=x^{\left(+\right)} \left(\cdot ;\, \vartheta ,\, \varepsilon \right)\right|_{\lambda =\varepsilon } $$\left(x^{\left(-\right)} \left(\cdot ;\tilde{\vartheta },\, \varepsilon \right)\right. $$\left. \left. =x^{\left(-\right)} \left(\cdot ;\, \vartheta ,\, \varepsilon \right)\right|_{\lambda =\varepsilon } \right)$, 
\\$\varepsilon \in \left(0,\, \bar{\varepsilon }\right]\bigcap \left(0,\, \tilde{\varepsilon }\right]\bigcap \left(0,\, 1\right)$, where $x^{\left(+\right)} \left(\cdot ;\vartheta ,\, \varepsilon \right)$ $\left(x^{\left(-\right)} \left(\cdot ;\vartheta ,\, \varepsilon \right)\right)$ was defined as a variation of the function $\bar{x}\left(\cdot \right)$ from \eqref{GrindEQ__2_4_} ((2.5)) allowing for \eqref{GrindEQ__2_1_} ((2.6)).

Considering the definition of the function $x^{\left(+\right)} \left(\cdot ;\tilde{\vartheta },\, \varepsilon \right)$ $\left(x^{\left(-\right)} \left(\cdot ;\tilde{\vartheta },\, \varepsilon \right)\right)$ and assuming $\xi =\eta ,$ similar to \eqref{GrindEQ__3_24_} ((3.25)) we have that for all $\varepsilon \in \left(0,\, \bar{\varepsilon }\right]\bigcap \left(0,\, \tilde{\varepsilon }\right]\bigcap \left(0,\, 2^{-1} \right]$  the following estimations are valid  
\begin{equation} \label{GrindEQ__3_27_} 
\max \left\{\left\| x^{\left(+\right)} \left(\cdot ;\, \tilde{\vartheta },\varepsilon \right)-\bar{x}\left(\cdot \right)\right\| _{C\left(I,\, R^{n} \right)} ,\, \left\| \dot{x}^{\left(+\right)} \left(\cdot ,\, \tilde{\vartheta },\varepsilon \right)-\dot{\bar{x}}\left(\cdot \right)\right\| _{L_{\infty } \left(I,\, R^{n} \right)} \, \, \right\}\le \left\| \eta \right\| _{R^{n} }  
\end{equation} 
\begin{equation} \label{GrindEQ__3_28_} 
\left(\max \left\{\left\| x^{\left(-\right)} \left(\cdot ;\, \tilde{\vartheta },\varepsilon \right)-\bar{x}\left(\cdot \right)\right\| _{C\left(I,\, R^{n} \right)} ,\, \left\| \dot{x}^{\left(-\right)} \left(\cdot ,\, \tilde{\vartheta },\varepsilon \right)-\dot{\bar{x}}\left(\cdot \right)\right\| _{L_{\infty } \left(I,\, R^{n} \right)} \, \, \right\}\le \left\| \eta \right\| _{R^{n} } \right).   
\end{equation} 
Note that here the inequality $\left|\left(\varepsilon -1\right)^{-1} \varepsilon \right|\le 1$ was considered for $\varepsilon \in \left(0,\, \, 2^{-1} \, \right]$.

Furthermore, carrying out similar reasoning stated in the proof of Theorem 3.4, allowing for estimations \eqref{GrindEQ__3_27_} and \eqref{GrindEQ__3_28_} the proof of Theorem 3.7 directly follows from Theorem 3.3. Theorem 3.7 is proved.

\section{Necessary conditions for a minimum in the presence of various degenerations on the interval.}

 In frequent cases, the Weierstrass condition and also the Legendre condition degenerate on some interval. Such a situation as an independent problem is studied in this section. It is important to note that the research of such cases allows to get analogues of statements \eqref{GrindEQ__3_6_}, \eqref{GrindEQ__3_14_} and \eqref{GrindEQ__3_22_} under significantly weakened assumptions on the smoothness of the integrant $L\left(\cdot \right)$ and the considered extremal of the problem \eqref{GrindEQ__1_1_}, \eqref{GrindEQ__1_2_}. Here as in section 3, the used research approach is based on the introduced special variations of the extremal of problem \eqref{GrindEQ__1_1_}, \eqref{GrindEQ__1_2_} and on increment formulas of functional \eqref{GrindEQ__1_1_} obtained in section 2.  

We prove the following theorems.

\textbf{ Theorem 4.1. }Let the integrant $L\left(\cdot \right)$ be continuously differentiable in totality of variables, furthermore, the admissible function $\bar{x}\left(\cdot \right)$ be an extremal in problem \eqref{GrindEQ__1_1_}, \eqref{GrindEQ__1_2_} and along it for the vectors $\eta \ne 0$ and $\left(\bar{\lambda }-1\right)^{-1} \bar{\lambda }\eta $, where $\bar{\lambda }\in \left(0,1\right)$, the Weierstrass condition degenerates at any point of the interval $\left(\bar{t}_{0} ,\, \bar{t}_{1} \right)\subset \left[t_{0} ,\, t_{1} \right]$, i.e. the following equalities hold: 
\begin{equation} \label{GrindEQ__4_1_} 
E\left(\bar{L}\right)\left(t,\, \eta \right)=E\left(\bar{L}\right)\, \left(t,\, \left(\bar{\lambda }-1\right)^{-1} \bar{\lambda }\eta \, \right)=0,\, \, \, \forall t\in \left(\bar{t}_{0} ,\, \bar{t}_{1} \right),   
\end{equation} 
where $E\left(\bar{L}\right)\left(\cdot \right)$ is determined from \eqref{GrindEQ__1_6_} and the interval $\left(\bar{t}_{0} ,\, \bar{t}_{1} \right)$ does not contain any angular point of function $\bar{x}\left(\cdot \right)$. Then: 

(i) if the extremal $\bar{x}\left(\cdot \right)$ is a strong local minimum in problem \eqref{GrindEQ__1_1_}, \eqref{GrindEQ__1_2_}, then the following equality is fulfilled:
\begin{equation} \label{GrindEQ__4_2_} 
\bar{\lambda }\Delta \bar{L}_{x}^{T} \left(t,\, \eta \right)\eta +\left(1-\bar{\lambda }\right)\Delta \bar{L}_{x}^{T} \left(t,\, \left(\bar{\lambda }-1\right)^{-1} \bar{\lambda }\eta \right)\eta =0,\, \, \forall t\in \left(\bar{t}_{0} ,\, \bar{t}_{1} \right),   
\end{equation} 
where $\Delta \bar{L}_{x} \left(\cdot \right)$  is determined by \eqref{GrindEQ__2_10_}; 

(ii) if the extremal $\bar{x}\left(\cdot \right)$ is a weak local minimum in problem \eqref{GrindEQ__1_1_}, \eqref{GrindEQ__1_2_}, then there exists a number $\delta >0$ such that for every point $\left(\eta ,\left(\bar{\lambda }-1\right)^{-1} \bar{\lambda }\eta ,\bar{\lambda }\right)\in B_{\delta } \left(0\right)\times B_{\delta } \left(0\right)\times \left(0,1\right)$ satisfying condition \eqref{GrindEQ__4_1_},  equality \eqref{GrindEQ__4_2_} is fulfilled.

 \textbf{Proof. }Since the admissible function $\bar{x}\left(\cdot \right)$ is an extremal in problem \eqref{GrindEQ__1_1_}, \eqref{GrindEQ__1_2_}, then along the function $\bar{x}\left(\cdot \right)$ formulas \eqref{GrindEQ__2_17_}-\eqref{GrindEQ__2_19_} are valid for the increment $\Delta _{\varepsilon }^{\left(+\right)} J\left(\bar{x}\left(\cdot \right),\vartheta \right)$, where $\vartheta =\left(\theta ,\, \lambda ,\, \xi \right)\in \left[t_{0} ,\, t_{1} \right)\times \left(0,\, 1\right)\times R^{n} \backslash \left\{0\right\}$ is an arbitrary fixed point, and $\varepsilon \in \left(0,\, \bar{\varepsilon }\right)$. Assume $\theta \in \left(\bar{t}_{0} ,\, \bar{t}_{1} \right)$, $\lambda =\bar{\lambda }$ and $\xi =\eta $, more exactly, $\vartheta =\bar{\vartheta }:=\left(\theta ,\, \bar{\lambda },\eta \right)\in \left(\bar{t}_{0} ,\, \bar{t}_{1} \right)\times \left(0,\, 1\right)\times R^{n} \backslash \left\{0\right\}$ and $\bar{\varepsilon }=\hat{\varepsilon }$, moreover $\hat{\varepsilon }<\bar{t}_{1} -\theta $. Then by virtue of \eqref{GrindEQ__1_6_}, \eqref{GrindEQ__2_1_} and \eqref{GrindEQ__2_2_}, the increment $\Delta _{\varepsilon }^{\left(+\right)} J\left(\bar{x}\left(\cdot \right),\bar{\vartheta }\right)$ takes the form 
\[\Delta _{\varepsilon }^{\left(+\right)} J\left(\bar{x}\left(\cdot \right),\bar{\vartheta }\right)=\int _{\theta }^{\theta +\bar{\lambda }\varepsilon }E\left(\bar{L}\right)\left(t,\, \eta \right)dt+\int _{\theta +\bar{\lambda }\varepsilon }^{\theta +\varepsilon }E\left(\bar{L}\right)\, \left(t,\, \left(\bar{\lambda }-1\right)^{-1} \bar{\lambda }\eta \right)dt+  \] 
\begin{equation} \label{GrindEQ__4_3_} 
+\hat{J}_{1}^{\left(+\right)} \left(\varepsilon ,\, \bar{\vartheta }\right)+\hat{J}_{2}^{\left(+\right)} \left(\varepsilon ,\, \bar{\vartheta }\right),\, \, \, \varepsilon \in \left(0,\, \hat{\varepsilon }\right].    
\end{equation} 
Here 
\[\hat{J}_{1}^{\left(+\right)} \left(\varepsilon ,\, \bar{\vartheta }\right)=\int _{\theta }^{\theta +\bar{\lambda }\varepsilon }\, \, \, \left[L\left(t,\bar{x}\left(t\right)+\left(t-\theta \right)\eta ,\, \dot{\bar{x}}\left(t\right)+\eta \right)-L\, \left(t,\bar{x}\left(t\right),\, \dot{\bar{x}}\left(t\right)+\eta \right)-\left(t-\theta \right)\, \bar{L}_{x}^{T} \left(t\right)\eta \right]dt ,\] 
\[\hat{J}_{2}^{\left(+\right)} \left(\varepsilon ,\, \bar{\vartheta }\right)=\int _{\theta +\bar{\lambda }\varepsilon }^{\theta +\varepsilon }\left[L\left(t,\bar{x}\left(t\right)+\frac{\bar{\lambda }}{\bar{\lambda }-1} \left(t-\theta -\varepsilon \right)\eta ,\, \dot{\bar{x}}\left(t\right)+\frac{\bar{\lambda }}{\bar{\lambda }-1} \eta \right)-L\, \left(t,\bar{x}\left(t\right),\, \dot{\bar{x}}\left(t\right)+\frac{\bar{\lambda }}{\bar{\lambda }-1} \eta \right)-\right. \]

 $\left. -\frac{\bar{\lambda }}{\bar{\lambda }-1} \left(t-\theta -\varepsilon \right)\, \bar{L}_{x}^{T} \left(t\right)\eta \right]\, dt.  $

Applying the Taylor formula and taking into account \eqref{GrindEQ__2_10_}, for $\hat{J}_{1}^{\left(+\right)} \left(\varepsilon ,\, \bar{\vartheta }\right)$ and$\hat{J}_{2}^{\left(+\right)} \left(\varepsilon ,\, \bar{\vartheta }\right)$  we have 
\begin{equation} \label{GrindEQ__4_4_} 
\hat{J}_{1}^{\left(+\right)} \left(\varepsilon ,\, \bar{\vartheta }\right)=\frac{1}{2} \varepsilon ^{2} \bar{\lambda }^{2} \Delta \bar{L}_{x}^{T} \left(\theta ,\eta \right)\eta +o \left(\varepsilon ^{2} \right),    
\end{equation} 
\begin{equation} \label{GrindEQ__4_5_} 
\hat{J}_{2}^{\left(+\right)} \left(\varepsilon ,\, \bar{\vartheta }\right)=\frac{1}{2} \varepsilon ^{2} \bar{\lambda }\left(1-\bar{\lambda }\right)\Delta \bar{L}_{x}^{T} \left(\theta ,\, \left(\bar{\lambda }-1\right)^{-1} \bar{\lambda }\eta \right)\eta +o \left(\varepsilon ^{2} \right).  
\end{equation} 
Further, since $\theta \in \left(\bar{t}_{0} ,\, \bar{t}_{1} \right)$, $\bar{\lambda }\in \left(0,\, 1\right)$ and $\hat{\varepsilon }<\bar{t}_{1} -\theta $, then for all $\varepsilon \in \left(0,\, \hat{\varepsilon }\right]$  the inclusions $\left[\, \theta ,\, \theta +\bar{\lambda }\varepsilon \, \right]\subset \left(\bar{t}_{0} ,\, \bar{t}_{1} \right)$ and $\left[\, \theta +\bar{\lambda }\varepsilon ,\, \theta +\varepsilon \, \right]\subset \left(\bar{t}_{0} ,\, \bar{t}_{1} \right)$ hold.  Therefore, by virtue of assumption \eqref{GrindEQ__4_1_} the first two terms in \eqref{GrindEQ__4_3_} vanish. Considering this and also \eqref{GrindEQ__4_4_} and \eqref{GrindEQ__4_5_}, from \eqref{GrindEQ__4_3_} we get 
\begin{equation} \label{GrindEQ__4_6_} 
\Delta _{\varepsilon }^{\left(+\right)} J\left(\bar{x}\left(\cdot \right),\, \bar{\vartheta }\right)=\frac{1}{2} \varepsilon ^{2} \left[\bar{\lambda }^{\, 2} \Delta L_{x}^{T} \left(\theta ,\, \eta \right)\eta +\bar{\lambda }\left(1-\bar{\lambda }\right)\Delta L_{x}^{T} \left(\theta ,\, \left(\bar{\lambda }-1\right)^{-1} \bar{\lambda }\eta \right)\, \eta \right]
+o\left(\varepsilon ^{2} \right),\, \, \varepsilon \in \left(0,\, \hat{\varepsilon }\right].      
\end{equation} 
Now, quite similarly, using \eqref{GrindEQ__2_6_}, \eqref{GrindEQ__2_7_}, \eqref{GrindEQ__2_10_} and \eqref{GrindEQ__2_25_}-\eqref{GrindEQ__2_27_} allowing for $\vartheta =\bar{\vartheta }=\left(\theta ,\, \bar{\lambda },\, \eta \right)$ and choosing the number $\tilde{\varepsilon }=\varepsilon ^{*} <\theta -\bar{t}_{0} $, we prove that by virtue of assumption \eqref{GrindEQ__4_1_} the increment$\Delta _{\varepsilon }^{\left(-\right)} J\left(\bar{x}\left(\cdot \right),\, \bar{\vartheta }\right)$ has the form 
\[\Delta _{\varepsilon }^{\left(-\right)} J\left(\bar{x}\left(\cdot \right),\, \vartheta \right)=-\frac{1}{2} \varepsilon ^{2} \left[\bar{\lambda }^{\, 2} \Delta L_{x}^{T} \left(\theta ,\, \eta \right)\eta +\bar{\lambda }\left(1-\bar{\lambda }\right)\Delta L_{x}^{T} \left(\theta ,\, \left(\bar{\lambda }-1\right)^{-1} \bar{\lambda }\eta \right)\eta \right]+\] 
\begin{equation} \label{GrindEQ__4_7_} 
+o\left(\varepsilon ^{2} \right),\, \, \varepsilon \in \left(0,\, \varepsilon ^{*} \right].      
\end{equation} 
Since the admissible function $\bar{x}(\cdot)$ is a strong local minimum in problem \eqref{GrindEQ__1_1_}, \eqref{GrindEQ__1_2_}, then $\Delta _{\varepsilon }^{\left(+\right)} J\left(\bar{x}\left(\cdot \right),\, \bar{\vartheta }\right)\ge 0$ and $\Delta _{\varepsilon }^{\left(-\right)} J\left(\bar{x}\left(\cdot \right),\, \bar{\vartheta }\right)\ge 0$ for all $\varepsilon \in \left(0,\, \hat{\varepsilon }\right]\bigcap \, \left(\, 0,\, \varepsilon ^{*} \, \right]$. Considering the last inequalities and arbitrariness of the point $\theta \in \left(\bar{t}_{0} ,\, \bar{t}_{1} \right)$, by virtue of  \eqref{GrindEQ__4_6_} and \eqref{GrindEQ__4_7_} we get the validity of the sought-for equality \eqref{GrindEQ__4_2_}, i.e. part (i) of Theorem 4.1 is proved. 

We now prove part (ii) of Theorem 4.1. Let us consider formulas \eqref{GrindEQ__4_6_} and \eqref{GrindEQ__4_7_} obtained for the increments $\Delta _{\varepsilon }^{\left(+\right)} J\left(\bar{x}\left(\cdot \right),\, \bar{\vartheta }\right)$ and $\Delta _{\varepsilon }^{\left(-\right)} J\left(\bar{x}\left(\cdot \right),\, \bar{\vartheta }\right)$. Here, by definition $\Delta _{\varepsilon }^{\left(+\right)} J\left(\bar{x}\left(\cdot \right),\, \bar{\vartheta }\right)=J\left(x^{\left(+\right)} \left(\cdot ;\, \bar{\vartheta },\, \varepsilon \right)\right)-J\left(\bar{x}\left(\cdot \right)\right)$ and $\Delta _{\varepsilon }^{\left(-\right)} J\left(\bar{x}\left(\cdot \right),\, \bar{\vartheta }\right)=J\left(x^{\left(-\right)} \left(\cdot ;\, \bar{\vartheta },\, \varepsilon \right)\right)-J\left(\bar{x}\left(\cdot \right)\right)$, where $x^{\left(+\right)} \left(\cdot ;\, \bar{\vartheta },\, \varepsilon \right)$ is determined from \eqref{GrindEQ__2_4_} allowing for  \eqref{GrindEQ__2_1_}, while $x^{\left(-\right)} \left(\cdot ;\, \bar{\vartheta },\, \varepsilon \right)$ is determined from \eqref{GrindEQ__2_5_} allowing for \eqref{GrindEQ__2_6_} at $\vartheta =\bar{\vartheta }=\left(\theta ,\, \eta ,\, \bar{\lambda }\right)$, i.e. at $\xi =\eta ,\, \, \lambda =\bar{\lambda }$ and $\theta \in \left(\bar{t}_{0} ,\, \bar{t}_{1} \right)$. 

By virtue of definition $x^{\left(+\right)}(\cdot ;\, \bar{\vartheta },\, \varepsilon)$ and $x^{\left(-\right)}(\cdot ;\, \bar{\vartheta },\, \varepsilon)$ allowing for \eqref{GrindEQ__2_2_} and \eqref{GrindEQ__2_7_} the following estimations are valid: 
\begin{equation} \label{GrindEQ__4_8_} 
\left\| x^{\left(+\right)} \left(\cdot ;\, \bar{\vartheta },\, \varepsilon \right)-\bar{x}\left(\cdot \right)\right\| _{C\left(I,\, R^{n} \right)} =\left\| h^{\left(+\right)} \left(\cdot ;\, \bar{\vartheta },\, \varepsilon \right)\right\| _{C\left(I,\, R^{n} \right)} \le \left\| \eta \right\| ,\, \, \varepsilon \in \left(0,\, \hat{\varepsilon }\right]\bigcap \left(0,\, 1\right),   
\end{equation} 
\begin{equation} \label{GrindEQ__4_9_} 
\left\| \dot{x}^{\left(+\right)} \left(\cdot ;\, \bar{\vartheta },\, \varepsilon \right)-\dot{\bar{x}}\left(\cdot \right)\right\| _{L_{\infty } \left(I,\, R^{n} \right)} =\left\| \dot{h}^{\left(+\right)} \left(\cdot ;\, \bar{\vartheta },\, \varepsilon \right)\right\| _{{}_{L_{\infty } } \left(I,\, R^{n} \right)} \le \max \left\{\, \left\| \eta \right\| \, _{R^{n} } ,\, \frac{\bar{\lambda }}{1-\bar{\lambda }} \left\| \eta \right\| \, _{R^{n} } \right\} 
\end{equation} 

 Similar estimations are valid for  $\left\| x^{\left(-\right)} \left(\cdot ;\, \bar{\vartheta },\, \varepsilon \right)-\bar{x}\left(\cdot \right)\right\| _{C\left(I,\, R^{n} \right)} $ and $\left\| \dot{x}^{\left(-\right)} \left(\cdot ;\, \bar{\vartheta },\, \varepsilon \right)-\dot{\bar{x}}\left(\cdot \right)\right\| _{L_{\infty } \left(I,\, R^{n} \right)} $ as well. Let the function $\bar{x}\left(\cdot \right)$ be a local minimum in problem \eqref{GrindEQ__1_1_}, \eqref{GrindEQ__1_2_} with $\hat{\delta }$-neighborhood. Then, considering the last estimations we confirm that for every point $\left(\eta ,\, \left(\bar{\lambda }-1\right)^{-1} \bar{\lambda }\eta ,\bar{\lambda }\right)\in B_{\hat{\delta }} \left(0\right)\times B_{\hat{\delta }} \left(0\right)\times \left(0,\, 1\right)$, satisfying condition \eqref{GrindEQ__4_1_}, the inequalities $\Delta _{\varepsilon }^{\left(+\right)} J\left(\bar{x}\left(\cdot \right),\, \bar{\vartheta }\right)\ge 0$ and $\Delta _{\varepsilon }^{\left(-\right)} J\left(\bar{x}\left(\cdot \right),\, \bar{\vartheta }\right)\ge 0$ are fulfilled. Based on these inequalities, allowing for \eqref{GrindEQ__4_6_}, \eqref{GrindEQ__4_7_} and arbitrariness of $\theta \in \left(\bar{t}_{0} ,\, \bar{t}_{1} \right)$, and also choosing $\delta =\hat{\delta }$, we get the proof of part (ii) of Theorem 4.1. So, Theorem 4.1 is completely proved.

\textbf{Theorem 4.2. }Let the functions $L\left(\cdot \right)$ and $L_{x} \left(\cdot \right)$ be continuously differentiable in totality of variables, and the admissible function $\bar{x}\left(\cdot \right)$ be an extremal of problem \eqref{GrindEQ__1_1_}, \eqref{GrindEQ__1_2_}, along it for the vectors $\eta \ne 0$ and $\left(\bar{\lambda }-1\right)^{-1} \bar{\lambda }\eta $, where $\bar{\lambda }\in \left(0,\, 1\right)$, the Weierstrass condition degenerates at any point of the interval $\left(\bar{t}_{0} ,\, \bar{t}_{1} \right)\subset \left[t_{0} ,\, t_{1} \right]$, i.e. condition \eqref{GrindEQ__4_1_} holds. Furthermore, let the function $\bar{x}\left(\cdot \right)$ be twice continuously differentiable on the interval $\left(\bar{t}_{0} ,\, \bar{t}_{1} \right)$. Then:

 (i) if the extremal $\bar{x}\left(\cdot \right)$ is a strong local minimum in problem \eqref{GrindEQ__1_1_}, \eqref{GrindEQ__1_2_}, then the following inequality is fulfilled: 
\begin{equation} \label{GrindEQ__4_10_} 
\eta ^{T} \left(\bar{\lambda }\bar{L}_{xx} \left(t,\, \eta \right)+\left(1-\bar{\lambda }\right)\, \bar{L}_{xx} \left(t,\, \left(\bar{\lambda }-1\right)^{-1} \bar{\lambda }\eta \right)\right)\, \eta -\frac{d}{dt} \Delta \bar{L}_{x}^{T} \, \left(t,\, \eta \right)\eta \ge 0,\, \, \forall t\in \left(\bar{t}_{0} ,\, \bar{t}_{1} \right), 
\end{equation} 
where $\bar{L}_{xx} \left(t,\, \cdot \right)$ and $\Delta \bar{L}_{x} \left(t,\, \cdot \right)$ are determined by \eqref{GrindEQ__2_9_} and \eqref{GrindEQ__2_10_}, respectively;

  (ii) if the extremal $\bar{x}\left(\cdot \right)$ is a weak local minimum in problem \eqref{GrindEQ__1_1_}, \eqref{GrindEQ__1_2_}, then there exists a number $\delta >0$ such that for every point  $\left(\eta ,\, \left(\bar{\lambda }-1\right)^{-1} \bar{\lambda }\eta ,\, \bar{\lambda }\right)\in B_{\delta } \left(0\right)\times B_{\delta } \left(0\right)\times \left(0,\, 1\right)$ satisfying the condition \eqref{GrindEQ__4_1_}, the inequality \eqref{GrindEQ__4_10_} is fulfilled.

 \textbf{Proof. }Let us\textbf{ }prove part (i) of Theorem 4.2. Obviously Theorem 4.1 is valid subject to the condition of Theorem 4.2. Let us use the formula determined by \eqref{GrindEQ__4_3_} for the increment $\Delta _{\varepsilon }^{\left(+\right)} J\left(\bar{x}\left(\cdot \right),\, \bar{\vartheta }\right)$, where $\bar{\vartheta }=\left(\theta ,\, \bar{\lambda },\, \eta \right)\in \left(\bar{t}_{0} ,\, \bar{t}_{1} \right)\times \left(0,\, 1\right)\times R_{n} \backslash \left\{0\right\}$ and $\varepsilon \in \left(0,\, \hat{\varepsilon }\right]$. For finding the increment $\Delta _{\varepsilon }^{\left(+\right)} J\left(\bar{x}\left(\cdot \right),\, \bar{\vartheta }\right)$ it suffices to find the expansion of the integrals $\hat{J}_{1}^{\left(+\right)} \left(\varepsilon ,\, \bar{\vartheta }\right)$ and $\hat{J}_{2}^{\left(+\right)} \left(\varepsilon ,\, \bar{\vartheta }\right)$ determined by \eqref{GrindEQ__4_4_} and \eqref{GrindEQ__4_5_}, with accuracy $o(\varepsilon ^{3})$. These integrals are calculated similar to \eqref{GrindEQ__2_22_} and \eqref{GrindEQ__2_23_}. More exactly, from \eqref{GrindEQ__4_4_} and \eqref{GrindEQ__4_5_}, using the Taylor formula, allowing for \eqref{GrindEQ__2_9_} and \eqref{GrindEQ__2_10_}, we get

\begin{equation} \label{GrindEQ__4_11_} 
\hat{J}_{1}^{\left(+\right)} \left(\varepsilon ,\, \bar{\vartheta }\right)=\frac{1}{2} \varepsilon ^{2} \bar{\lambda }^{2} \Delta \bar{L}_{x}^{T} \left(\theta ,\, \eta \right)\eta +\frac{1}{6} 
\varepsilon^{2} [2 \bar{\lambda }^{3} \frac{d}{dt} \Delta \bar{L}_{x}^{T} \left(\theta ,\, \eta \right)\eta +\bar{\lambda }^{3} \eta ^{T} \bar{L}_{xx} \left(\theta ,\, \eta \right)\eta 
]+o\left(\varepsilon ^{3} \right),      
\end{equation} 
\[\hat{J}_{2}^{\left(+\right)} \left(\varepsilon ,\, \bar{\vartheta }\right)=\frac{1}{2} \varepsilon ^{2} \bar{\lambda }\left(1-\bar{\lambda }\right)\Delta \bar{L}_{x}^{T} \left(\theta ,\frac{\bar{\lambda }}{\bar{\lambda }-1} \, \eta \right)\eta +\] 
\begin{equation} \label{GrindEQ__4_12_} 
\left. +\frac{1}{6} \varepsilon ^{2} \left[\bar{\lambda }\left(1-\bar{\lambda }\right)\left(1+2\bar{\lambda }\right)\right. \frac{d}{dt} \Delta \bar{L}_{x}^{T} \left(\theta ,\frac{\bar{\lambda }}{\bar{\lambda }-1} \, \eta \right)\eta +\bar{\lambda }^{2} \left(1-\bar{\lambda }\right)\eta ^{T} \bar{L}_{xx} \left(\theta ,\frac{\bar{\lambda }}{\bar{\lambda }-1} \, \eta \right)\eta \right]+o\left(\varepsilon ^{3} \right).    
\end{equation} 

Since $\theta \in \left(t_{0} ,\, \bar{t}_{1} \right)$, $\bar{\lambda }\in \left(0,\, 1\right)$, $\hat{\varepsilon }<t_{1} -\theta $ and $\bar{x}\left(\cdot \right)$ is a strong local minimum in problem \eqref{GrindEQ__1_1_}, \eqref{GrindEQ__1_2_}, firstly, by assumption \eqref{GrindEQ__4_1_} for all $\varepsilon \in \left(0,\, \hat{\varepsilon }\right]$ the first two terms in \eqref{GrindEQ__4_3_}  vanish; secondly, for all $\varepsilon \in \left(0,\, \hat{\varepsilon }\right]$ the increment $\Delta _{\varepsilon }^{\left(+\right)} J\left(\bar{x}\left(\cdot \right),\, \bar{\vartheta }\right)$ is non-negative, and by virtue of statement \eqref{GrindEQ__4_2_} of Theorem 4.1 for all $t\in \left(\bar{t}_{0} ,\, \bar{t}_{1} \right)$ the equality $\Delta \bar{L}_{x}^{T} \left(t,\, \left(\bar{\lambda }-1\right)^{-1} \bar{\lambda }\eta \right)\, \eta =\left(\bar{\lambda }-1\right)^{-1} \bar{\lambda }\Delta \bar{L}_{x}^{T} \left(t,\, \eta \right)\eta $ is valid. By means of \eqref{GrindEQ__4_11_}, \eqref{GrindEQ__4_12_} and the last statements, allowing for arbitrariness of $\theta \in \left(\bar{t}_{0} ,\, \bar{t}_{1} \right)$, the sought-for inequality \eqref{GrindEQ__4_10_} follows from \eqref{GrindEQ__4_3_} i.e. part (i) of Theorem 4.2 is proved. 

The proof of part (ii) of Theorem 4.2, allowing for estimations \eqref{GrindEQ__4_8_}, \eqref{GrindEQ__4_9_} and definition of a weak local minimum of the function $\bar{x}\left(\cdot \right)$, follows from part (i) of Theorem 4.2. By the same token, Theorem 4.2 is completely proved. 

Finally, we prove the following theorem. 

\textbf{Theorem 4.3. }Let the integrant $L\left(\cdot \right)$ be continuously differentiable in totality of variables, and partial derivatives of the form $L_{xx} \left(\cdot \right),\, L_{x\dot{x}} \left(\cdot \right),\, L_{\dot{x}\dot{x}} \left(\cdot \right)$ and $L_{\dot{x}\dot{x}\dot{x}} \left(\cdot \right)$ be continuous in totality of variables. Furthermore, let the admissible function $\bar{x}\left(\cdot \right)$ be an extremal of problem \eqref{GrindEQ__1_1_}, \eqref{GrindEQ__1_2_}, and along it for the vector $\eta \ne 0$ the Weierstrass and Legendre conditions degenerate at any point of the interval $\left(\bar{t}_{0} ,\, \bar{t}_{1} \right)\subset \left[t_{0} ,\, t_{1} \right]$, i.e. the following equalities hold
\begin{equation} \label{GrindEQ__4_13_} 
E\left(\bar{L}\right)\left(t,\, \eta \right)=\eta ^{T} \bar{L}_{\dot{x}\dot{x}} \left(t\right)\eta =0,\, \, \forall t\in \left(\bar{t}_{0} ,\, \bar{t}_{1} \right),     
\end{equation} 
where the interval $\left(\bar{t}_{0} ,\, \bar{t}_{1} \right)$ does not contain any angular point of extremal $\bar{x}\left(\cdot \right)$. Then:

 (i) if the extremal $\bar{x}\left(\cdot \right)$ is a strong local minimum in problem \eqref{GrindEQ__1_1_}, \eqref{GrindEQ__1_2_}, then the following equality is fulfilled
\begin{equation} \label{GrindEQ__4_14_} 
\eta ^{T} \left(L_{x} \left(t,\, \bar{x}\left(t\right),\, \dot{\bar{x}}\left(t\right)+\eta \right)-\bar{L}_{x} \left(t\right)-\bar{L}_{x\dot{x}} \left(t\right)\, \eta \right)=0,\, \, \forall t\in \left(\bar{t}_{0} ,\, \, \bar{t}_{1} \right);    
\end{equation} 

(ii) if the extremal $\bar{x}\left(\cdot \right)$ is a weak local minimum in problem \eqref{GrindEQ__1_1_}, \eqref{GrindEQ__1_2_}, the there exists a number $\delta >0$ such that for every point $\eta \in B_{\delta } \left(0\right)$ satisfying condition \eqref{GrindEQ__4_13_}, the equality \eqref{GrindEQ__4_14_} is fulfilled.

\textbf{Proof. }Prove part (i) of Theorem 4.3. Use the formulas \eqref{GrindEQ__2_39_}-\eqref{GrindEQ__2_41_}, obtained for the increment $\Delta _{\varepsilon }^{\left(+\right)} J\left(\bar{x}\left(\cdot \right),\, \tilde{\vartheta }\right)=J_{1}^{\left(+\right)} \left(\varepsilon ,\, \tilde{\vartheta }\right)+J_{2}^{\left(+\right)} \left(\varepsilon ,\, \tilde{\vartheta }\right),\, \varepsilon \in \left(0,\, \bar{\varepsilon }\right]\bigcap \left(0,\, 1\right)$. Assume $\tilde{\vartheta }=\hat{\vartheta }$. Here $\tilde{\vartheta }=\hat{\vartheta }:=\left(\theta ,\, \varepsilon ,\, \eta \right)\in \left(\bar{t}_{0} ,\, \bar{t}_{1} \right)\times \left(0,\, \hat{\varepsilon }\right]\times R^{n} \backslash \left\{0\right\}$, where $\hat{\varepsilon }=\min \left\{\bar{\varepsilon },\, \bar{t}_{1} -\theta ,\, 1\right\}$. Then by the Taylor formula, allowing for $\hat{\vartheta }$ and $\hat{\varepsilon }$, firstly, from \eqref{GrindEQ__2_40_} for  $J_1^{\left(+\right)} \left(\varepsilon ,\, \hat{\vartheta }\right)$, $\varepsilon \in \left(0,\, \hat{\varepsilon }\right]$, we get 
\[J_{1}^{\left(+\right)} \left(\varepsilon ,\, \hat{\vartheta }\right)=\int _{\theta }^{\theta +\varepsilon ^{2} }E\left(\bar{L}\right)\left(t,\, \eta \right)dt+\frac{\varepsilon ^{2} }{2\left(\varepsilon -1\right)^{2} } \int _{\theta +\varepsilon ^{2} }^{\theta +\varepsilon }\eta ^{T} \bar{L}_{\dot{x}\dot{x}} \left(t\right)\, \eta \, dt-  \] 
\begin{equation} \label{GrindEQ__4_15_} 
-\frac{\varepsilon ^{4} }{6\left(1-\varepsilon \right)^{2} } \left(\eta ^{T} \bar{L}_{\dot{x}\dot{x}} \left(\theta \right)\eta \right)_{\dot{x}}^{T} \eta +o\left(\varepsilon ^{4} \right);     
\end{equation} 

secondly, taking into account \eqref{GrindEQ__2_41_}, similar to \eqref{GrindEQ__2_43_} for $J_{2}^{\left(+\right)} \left(\varepsilon ,\, \hat{\vartheta }\right)$, $\varepsilon \in \left(0,\, \hat{\varepsilon }\right]$ we have 
\begin{equation} \label{GrindEQ__4_16_} 
J_{2}^{\left(+\right)} \left(\varepsilon ,\, \hat{\vartheta }\right)=\frac{1}{2} \varepsilon ^{4} \left[\eta ^{T} \, \left(\bar{L}_{x} \left(\theta ,\, \eta \right)-\bar{L}_{x} \left(\theta \right)\right)-\eta ^{T} \bar{L}_{x\dot{x}} \left(\theta \right)\eta \right]+o\left(\varepsilon ^{4} \right),\, \forall \varepsilon \in \left(0,\, \hat{\varepsilon }\right].   
\end{equation} 

Further, since $\bar{x}\left(\cdot \right)$ is a strong local minimum in problem \eqref{GrindEQ__1_1_}, \eqref{GrindEQ__1_2_} and \eqref{GrindEQ__4_13_} is fulfilled, then by Theorem 3.6, and also definition of the point $\theta $ and the number $\hat{\varepsilon }$ for the increment $\Delta _{\varepsilon }^{\left(+\right)} J\left(\bar{x}\left(\cdot \right),\, \hat{\vartheta }\right)$ allowing for \eqref{GrindEQ__4_15_} and \eqref{GrindEQ__4_16_} the following inequality is valid: 
\begin{equation} \label{GrindEQ__4_17_} 
\Delta _{\varepsilon }^{\left(+\right)} J\left(\bar{x}\left(\cdot \right),\, \hat{\vartheta }\right)=\frac{1}{2} \varepsilon ^{4} \left[\eta ^{T} \, \left(\bar{L}_{x} \, \left(\theta ,\, \eta \right)-\bar{L}_{x} \left(\theta \right)\right)-\eta ^{T} \bar{L}_{x\dot{x}} \left(\theta \right)\eta \right]+
o\left(\varepsilon ^{4} \right)\ge 0,\, \, \forall \varepsilon \in \left(0,\, \hat{\varepsilon }\right].       
\end{equation} 

Quite similarly, using \eqref{GrindEQ__2_44_}-\eqref{GrindEQ__2_46_} allowing for
$\tilde{\vartheta }=\hat{\vartheta }:=\left(\theta ,\, \varepsilon ,\, \eta \right)\in \left(\bar{t}_{0} ,\, \bar{t}_{1} \right)\times \left(0,\, \varepsilon ^{*} \right]\times R^{n} \backslash \left\{0\right\}$, where $\varepsilon ^{*} =\left\{\tilde{\varepsilon },\, \theta -\bar{t}_{0} ,\, 1\right\}$, is it easy to show that for the increment $\Delta _{\varepsilon }^{\left(-\right)} J\left(\bar{x}\left(\cdot \right),\, \hat{\vartheta }\right)$ the inequality of the form 
\[\Delta _{\varepsilon }^{\left(-\right)} J\left(\bar{x}\left(\cdot \right),\, \hat{\vartheta }\right)=-\frac{1}{2} \varepsilon ^{4} \left[\eta ^{T} \, \left(\bar{L}_{x} \, \left(\theta ,\, \eta \right)-\bar{L}_{x} \left(\theta \right)\right)-\eta ^{T} \bar{L}_{x\dot{x}} \left(\theta \right)\eta \right]+\] 
\begin{equation} \label{GrindEQ__4_18_} 
o\left(\varepsilon ^{4} \right)\ge 0,\, \, \forall \varepsilon \in \left(0,\, \varepsilon ^{*} \right)     
\end{equation} 

is valid.

From inequalities \eqref{GrindEQ__4_17_} and \eqref{GrindEQ__4_18_}, allowing for arbitrariness of the point $\theta \in \left(\bar{t}_{0} ,\, \bar{t}_{1} \right)$, the sought-for equality \eqref{GrindEQ__4_14_} follows, i.e. part (i) of Theorem 4.3 is proved.

Prove part (ii) of Theorem 4.3.  Consider the increments  $\Delta _{\varepsilon }^{\left(+\right)} J\left(\bar{x}\left(\cdot \right),\, \hat{\vartheta }\right)=J\left(x^{\left(+\right)} \left(\cdot ,\, \hat{\vartheta },\, \varepsilon \right)\right)-J\left(\bar{x}\left(\cdot \right)\right)$ and $\Delta _{\varepsilon }^{\left(-\right)} J\left(\bar{x}\left(\cdot \right),\, \hat{\vartheta }\right)=J\left(x^{\left(-\right)} \left(\cdot ,\, \hat{\vartheta },\, \varepsilon \right)\right)-J\left(\bar{x}\left(\cdot \right)\right)$, determined above while proving part (i) of Theorem 4.3. Here by virtue of \eqref{GrindEQ__2_34_} and \eqref{GrindEQ__2_35_} allowing for $\hat{\vartheta} =\left(\theta ,\, \varepsilon ,\, \eta \right)$ we have $x^{\left(+\right)} \left(\cdot ,\, \hat{\vartheta },\, \varepsilon \right)=\left. x^{\left(+\right)} \left(\cdot ;\left(\theta ,\, \lambda ,\, \eta \right),\varepsilon \right)\right|_{\lambda =\varepsilon } $ and $x^{\left(-\right)} \left(\cdot ,\, \hat{\vartheta },\, \varepsilon \right)=\left. x^{\left(-\right)} \left(\cdot ;\left(\theta ,\, \lambda ,\, \eta \right),\varepsilon \right)\right|_{\lambda =\varepsilon } $, where $\varepsilon \in \left(0,\, \hat{\varepsilon }\right),\, $$x^{\left(+\right)} \left(\cdot ;\left(\theta ,\, \lambda ,\, \eta \right),\varepsilon \right)$ and $x^{\left(-\right)} \left(\cdot ;\left(\theta ,\, \lambda ,\, \eta \right),\varepsilon \right)$ are determined from \eqref{GrindEQ__2_4_} and \eqref{GrindEQ__2_5_} allowing for \eqref{GrindEQ__2_1_}, \eqref{GrindEQ__2_6_}, $\lambda =\varepsilon$ and $\xi =\eta $. Therefore, by virtue of estimations \eqref{GrindEQ__4_8_} and \eqref{GrindEQ__4_9_} for all $\varepsilon \in \left(0,\, 2^{-1} \, \right]\bigcap \left(0,\, \hat{\varepsilon }\right)$ the following estimations are valid: 
\begin{equation} \label{GrindEQ__4_19_} 
\max \left\{\left\| x^{\left(+\right)} \left(\cdot ;\, \hat{\vartheta },\, \varepsilon \right)-\bar{x}\left(\cdot \right)\right\| _{C\left(I,\, R^{n} \right)} ,\, \left\| \dot{x}^{\left(+\right)} \left(\cdot ;\, \hat{\vartheta },\, \varepsilon \right)-\dot{\bar{x}}\left(\cdot \right)\right\| _{L_{\infty } \left(I,\, R^{n} \right)} \right\}\le \left\| \eta \right\| _{R^{n} }  
\end{equation} 
In a similar way we have estimations of the form
\begin{equation} \label{GrindEQ__4_20_} 
\max \left\{\left\| x^{\left(-\right)} \left(\cdot ;\, \hat{\vartheta },\, \varepsilon \right)-\bar{x}\left(\cdot \right)\right\| _{C\left(I,\, R^{n} \right)} ,\, \left\| \dot{x}^{\left(-\right)} \left(\cdot ;\, \hat{\vartheta },\, \varepsilon \right)-\dot{\bar{x}}\left(\cdot \right)\right\| _{L_{\infty } \left(I,\, R^{n} \right)} \right\}\le \left\| \eta \right\| _{R^{n} } .   
\end{equation} 

Let the function $\bar{x}\left(\cdot \right)$ be a weak local minimum in problem \eqref{GrindEQ__1_1_}, \eqref{GrindEQ__1_2_} with $\hat{\delta }$-neighborhood. Then, considering estimation \eqref{GrindEQ__4_19_} and \eqref{GrindEQ__4_20_}, we confirm that for every point $\eta \in B_{\hat{\delta }} \left(0\right)$ satisfying the condition \eqref{GrindEQ__4_13_}, and for all  $\varepsilon \in \left(0,\, \, 2^{-1} \right]\bigcap \left(0,\, \hat{\varepsilon }\right)\bigcap \left(0,\, \varepsilon ^{*} \right)$ the inequalities $\Delta _{\varepsilon }^{\left(+\right)} J\left(\bar{x}\left(\cdot \right),\, \hat{\vartheta }\right)\ge 0$, $\Delta _{\varepsilon }^{\left(-\right)} J\left(\bar{x}\left(\cdot \right),\, \hat{\vartheta }\right)\ge 0$ are fulfilled. Based on these inequalities, allowing for \eqref{GrindEQ__4_17_}, \eqref{GrindEQ__4_18_} and arbitrariness of $\theta =\left(\bar{t}_{0} ,\, \bar{t}_{1} \right)$, and also choosing $\delta =\hat{\delta }$, we get the proof of part (ii) of Theorem 4.3. Thus, Theorem 4.3 is completely proved. 

 \textbf{Remark 4.1. }We consider a classical variational problem, for example with a free right end. We call it problem (A).  If the admissible function $x^{*} \left(\cdot \right)$ is a strong (weak) local minimum in problem (A), then obviously, this function gives strong (weak) local minimum to the functional \eqref{GrindEQ__1_1_} with the boundary conditions $x\left(t_{0} \right)=x_{0} $, $x\left(t_{1} \right)=x^{*} :=x^{*} \left(t_{1} \right)$. Therefore, all the statements of Theorem 3.1-3.7 and 4.1-4.3 are necessary conditions for a minimum of problem (A) as well.

\section{Discussions and examples} 

We consider a particular problem, more exactly the following problem with a free right end of the form
\begin{equation} \label{GrindEQ__5_1_} 
J\left(x\left(\cdot \right)\right)=\int _{t_{0} }^{t_{1} }L^{*} \left(t,\, x\left(t\right),\, \dot{x}\left(t\right)\right)\, dt\to {\mathop{\min }\limits_{x\left(\cdot \right)}}  ,    
\end{equation} 
\begin{equation} \label{GrindEQ__5_2_} 
x\left(t_{0} \right)=x_{0} ,\, \, x\left(\cdot \right)\in PC^{1} \left(\left[t_{0} ,\, t_{1} \right],\, R^{n} \right).    
\end{equation} 

Here $L^{*} \left(t,\, x,\, \dot{x}\right)=\dot{x}^{T} A\left(t\right)\dot{x}+2\dot{x}^{T} B\left(t\right)x+x^{T}C\left(t\right)x$, where $A\left(\cdot \right),\, B\left(\cdot \right)$ and $C\left(\cdot \right)$ are $n\times n$ -continuously differentiable functions.

Obviously, for every admissible function $\bar{x}\left(\cdot \right)$ (i.e. the function satisfying condition (5.2)) and the vector $\eta \in R^{n} $ we have the equality $E\left(\bar{L}\right)\left(t,\, \eta \right)=\eta ^{T} \bar{L}_{\dot{x}\dot{x}} \left(t\right)\eta ,\, t\in \left[t_{0} ,\, t_{1} \right]$, where $E\left(\bar{L}\right)\left(t,\, \eta \right)$ is determined by \eqref{GrindEQ__1_6_}, and $\bar{L}_{\dot{x}\dot{x}} \left(t\right)\equiv A\left(t\right)$. Hence, we have that in problem \eqref{GrindEQ__5_1_}, \eqref{GrindEQ__5_2_} that from every degeneration of the forms \eqref{GrindEQ__3_1_}, \eqref{GrindEQ__3_2_}, \eqref{GrindEQ__3_5_}, \eqref{GrindEQ__3_17_}. \eqref{GrindEQ__3_18_}, \eqref{GrindEQ__3_21_}, \eqref{GrindEQ__4_1_} and \eqref{GrindEQ__4_13_} it follows the degeneration of the Legendre condition and vice versa. Taking into account this property and Remark 4.1 when solving problem \eqref{GrindEQ__5_1_}, \eqref{GrindEQ__5_2_}, it is easy to conclude that the Kelly condition (see [8, p. 111]) and condition \eqref{GrindEQ__4_10_} coincide, namely, we have: if the admissible function $\bar{x}\left(\cdot \right)$ is a strong local minimum in problem \eqref{GrindEQ__5_1_}, \eqref{GrindEQ__5_2_} and for the vector $\eta \ne 0$ the Legendre condition degenerates at any point of the interval $\left(\bar{t}_{0} ,\, \bar{t}_{1} \right)\subset \left[t_{0} ,\, t_{1} \right]$, i.e. $\eta ^{T} \, A\left(t\right)\eta =0$, $\forall t\in \left(\bar{t}_{0} ,\, \bar{t}_{1} \right)$, then the inequality $\eta ^{T} \left[C\left(t\right)-\dot{B}\left(t\right)\right]\, \eta \ge 0$, $\forall t\in \left(\bar{t}_{0} ,\, \bar{t}_{1} \right)$ is valid.

Further, unlike the necessary condition for a minimum \eqref{GrindEQ__4_10_}, the Kelly condition is obtained only by means of the second variation of the functional when solving a singular optimal control problem.

So, based on what has been said, allowing for   below given Example 5.1 we can say that necessary condition for a minimum \eqref{GrindEQ__4_10_} is one reinforced variant of the Kelly condition in problem \eqref{GrindEQ__1_1_}, \eqref{GrindEQ__1_2_}.

Carrying out similar reasonings, we can confirm that in problem \eqref{GrindEQ__1_1_}, \eqref{GrindEQ__1_2_} the minimum condition \eqref{GrindEQ__4_2_}, and also minimum condition \eqref{GrindEQ__4_14_} are reinforced variant of the known equality type necessary condition obtained in \cite{8}. 

 We also note that, following the below given Example 5.1, we come to the conclusion that necessary condition \eqref{GrindEQ__4_10_} is not a corollary of the result of the paper \cite{7}, and it is more constructive compared to \cite{7}, more exactly, no need to solve the matrix differential equation. The similar conclusion refers also to necessary conditions obtained in Section 3 when comparing with appropriate results of the monography \cite{8}. Simple comparison and below given Examples 5.2 and 5.3 show that Theorem 3.3, 3.6 and 3.7 are sharpening and refinement of appropriate results of the paper \cite{3}.

It is important to note that degeneration of the form \eqref{GrindEQ__3_21_}, generally speaking does not follow from degeneration of the form \eqref{GrindEQ__3_5_}, and vice versa (see Example 5.2 and 5.4). Therefore, based on this proposition and formulation of the proved theorems in Section 3 and 4, it is easy to confirm that every condition for a minimum obtained by us has its own application area.

 \textbf{Example 5.1. }We consider a problem with a free right end of the form
\begin{equation} \label{GrindEQ__5_3_} 
J\left(x\left(\cdot \right)\right)=\int _{0}^{1}x^{2}  \left(1-\dot{x}^{2} \right)dt\to {\mathop{\min }\limits_{x\left(\cdot \right)}} ,\, \, x\left(0\right)=0,    
\end{equation} 
where $x\left(\cdot \right)\subset KC^{1} \left(\, \left[0,\, 1\right],\, R\right)$, $L\left(t,\, x,\, \dot{x}\right)=x^{2} \left(1-\dot{x}^{2} \right)$.

 Consider the admissible function $\bar{x}\left(t\right)=0,\, t\in \left[0,\, 1\right]$ for a minimum. Along this function allowing for \eqref{GrindEQ__1_6_}, \eqref{GrindEQ__2_9_} and \eqref{GrindEQ__2_10_} for all $\xi \in R$ and $t\in \left[0,1\right]$ we have 
\[\bar{L}\left(t\right)=L\left(t,\, \bar{x}\left(t\right),\, \dot{\bar{x}}\left(t\right)+\xi \right)\equiv 0,\, \, \, \, \bar{L}_{x} \left(t\right)=\bar{L}_{\dot{x}} \left(t\right)=\Delta \bar{L}_{x} \left(t,\, \xi \right)=L_{x} \left(t,\, \bar{x}\left(t\right),\, \dot{\bar{x}}\left(t\right)+\xi \right)\equiv 0,\] 
\[\bar{L}_{xx} \left(t,\, \xi \right)=2 \left(1-\xi ^{2} \right).\] 
Hence it is seen that the function $\bar{x}\left(\cdot \right)=0$ is an extremal of problem \eqref{GrindEQ__5_3_} and satisfies the Weierstrass condition \eqref{GrindEQ__1_5_}, and also along it for all $\left(\bar{\lambda },\, \eta \right)\in \left(0,1\right)\times R$ at any points of the interval (0,1) the assumption \eqref{GrindEQ__4_1_} is fulfilled. Therefore, considering Remark 4.1, i.e. possibility of application of Theorem 4.2, we have:

(a) condition \eqref{GrindEQ__4_10_} takes the form: $2 \eta ^{2} \left(1-\frac{\bar{\lambda }}{1-\bar{\lambda }} \eta ^{2} \right)\ge 0$  for all $\left(\bar{\lambda },\, \eta \right)\in \left(0,1\right)\times R$, and therefore for example for $\left(\bar{\lambda },\, \eta \right)=\left(\frac{1}{2} ,\, 2\right)$ it is not fulfilled. It means that the admissible function $\bar{x}\left(\cdot \right)=0$ by virtue of necessary minimum condition \eqref{GrindEQ__4_10_} is not a strong local minimum in problem \eqref{GrindEQ__5_3_}; 

(b) obviously, statement of part (ii) of Theorem 4.2 is fulfilled for $\delta =1$, and therefore an extremal $\bar{x}\left(\cdot \right)=0$ can be a weak local minimum in problem \eqref{GrindEQ__5_3_}. 

 Naturally, it is important to compare the necessary condition for a minimum \eqref{GrindEQ__4_10_} with the known necessary conditions for optimality of singular controls, for example, with the Kelly condition and the optimality condition in the paper \cite{7}.

 For that we consider the following problem equivalent to problem \eqref{GrindEQ__5_3_}:
\begin{equation} \label{GrindEQ__5_4_} 
\left\{\begin{array}{c} {\dot{x}_{1} =u,\, \, \, \, \, \, \, \, \, \, \, \, \, \, \, \, \, \, \, \, } \\ {\dot{x}_{2} =x_{1}^{2} \left(1-u^{2} \right),\, } \end{array}\, \, \, \, \, \begin{array}{c} {\, \, u\in R,\, \, \, } \\ {t\in \left[0,\, 1\right]\, \, \, \, \, x_{1} \left(0\right)=x_{2} \left(0\right)=0,} \end{array}\right.   x_{2} \left(1\right)\to \min , 
\end{equation} 
where $x_{1} \left(\cdot \right)=x\left(\cdot \right)$ is a sought-for admissible function in problem \eqref{GrindEQ__5_3_}. We study the admissible process $\left(\bar{u}\left(\cdot \right),\, \left(\bar{x}_{1} \left(\cdot \right),\, \bar{x}_{2} \left(\cdot \right)\right)^{T} \right)=\left(0,\, \left(0,0\right)^{T} \right)$ for optimality in problem \eqref{GrindEQ__5_4_}. Along this process we have: 

$\bar{\psi }\left(\cdot \right)=\left(\bar{\psi }_{1} \left(\cdot \right),\bar{\psi }_{2} \left(\cdot \right)\right)^{T} =\left(0,\, -1\right)^{T} $, $H\left(\bar{\psi }\left(\cdot \right),\, x_{1} ,\, x_{2} ,\, u,\, t\right)=-x_{1}^{2} \left(1-u^{2} \right)$,   

\noindent $H_{uu} \left(\bar{\psi }\left(\cdot \right),\, \bar{x}_{1} \left(\cdot \right),\, \bar{x}_{2} \left(\cdot \right),\, \bar{u}\left(\cdot \right),\, t\right)=0$, $H\left(\bar{\psi }\left(\cdot \right),\, \bar{x}_{1} \left(\cdot \right),\, \bar{x}_{2} \left(\cdot \right),\vartheta ,\, t\right)=H\left(\bar{\psi }\left(\cdot \right),\, \bar{x}_{1} \left(\cdot \right),\, \bar{x}_{2} \left(\cdot \right),\, \bar{u}\left(\cdot \right),\, t\right)=0$ 

for all $\vartheta \in R$ and $t\in \left[0,1\right]$, where $H\left(\cdot \right)$ is a Hamilton-Pontryagin function \cite{31}. Hence we easily get that the control $\bar{u}\left(\cdot \right)=0$ is singular in the classical sense (see, for example, [8, p. 28]), and also singular in Pontryagin's sense (see [8, p. 26]). Carrying out simple calculations, along the singular control $\overline{u}\left(\cdot \right)=0$ the Kelly condition takes the form: $2\vartheta ^{2} \ge 0$, for all $\vartheta \in R$. Further, since the control $\bar{u}\left(\cdot \right)=0$ is also singular in the sense of Pontryagin, then using the necessary optimality condition in the paper \cite{7}, we arrive at the inequality $2\left(t-1\right)\vartheta ^{2} \le 0$ for all $\vartheta \in R$ and $t\in \left[0,\, 1\right]$.

 As can be seen, the Kelly condition and also the condition in the paper \cite{7} along the control $\bar{u}\left(\cdot \right)=0$ are fulfilled, and both of these necessary conditions for optimality keep the control $\bar{u}\left(\cdot \right)=0$ among the contenders for optimality. Consequently, in the considered case i.e. while studying the extremal $\bar{x}\left(\cdot \right)=0$ of problem \eqref{GrindEQ__5_3_} for a minimum, each of  these necessary conditions for optimality is weak compared to minimum condition \eqref{GrindEQ__4_10_}.

\textbf{Example 5.2. }Consider the problem 
\begin{equation} \label{GrindEQ__5_5_} 
\int _{0}^{1}\left(\left(\dot{x}_{1} -\dot{x}_{2}^{2} \right)^{4} +x_{1} \dot{x}_{2}^{2} \right)dt\to {\mathop{\min }\limits_{x\left(\cdot \right)}} ,\, \, x_{i} \left(0\right)=x_{i} \left(1\right) =0,\, \, i=1,\, 2,    
\end{equation} 
where $x=\left(x_{1} ,\, x_{2} \right)^{T} ,\, \, L\left(\cdot \right)=\left(\dot{x}_{1} -\dot{x}_{2} \right)^{4} +x_{1} \dot{x}_{2}^{2} $.

 Obviously, the admissible function $\bar{x}\left(\cdot \right)=0$ is an extremal in problem \eqref{GrindEQ__5_5_} and along it the degeneration in the form \eqref{GrindEQ__4_13_} is fulfilled for all $t\in \left[0,\, 1\right]$  and for any vectors $\eta =\left(\eta _{1} ,\, \eta _{2} \right)^{T} $ such that $\eta _{1} =\eta _{2}^{2} $, $\eta _{2} \in R$,  since $E\left(\bar{L}\right)\left(t,\, \eta \right)=\left(\eta _{1} -\eta _{2}^{2} \right)^{4} $ and $\bar{L}_{\dot{x}\dot{x}}\left(t\right)=0,\, \, t\in \left[0,\, 1\right]$. Also we have: 
\begin{equation} \label{GrindEQ__5_6_} 
\eta ^{T} \left(L_{x} \left(t,\, \bar{x}\left(t\right),\, \dot{\bar{x}}\left(t\right)+\eta \right)-\bar{L}_{x} \left(t\right)-\bar{L}_{x\dot{x}} \left(t\right)\eta \right)=\eta _{2}^{4} ,    
\end{equation} 
where $\eta =\left(\eta _{2}^{2} ,\, \eta _{2} \right)^{T} ,\, \, \eta _{2} \in R$.

 Hence we get that Theorems 2.1 and 2.3 of the paper \cite{3} keep the extremal $\bar{x}\left(\cdot \right)=0$ among the contenders for a minimum and even for a weak local minimum in problem \eqref{GrindEQ__5_5_}. However, from \eqref{GrindEQ__5_6_} it can be seen that necessary condition for a minimum \eqref{GrindEQ__4_14_} is not fulfilled for all $\eta_2 \ne 0$, so, by virtue of Theorem 4.3 we have that the extremal $\bar{x}\left(\cdot \right)=0$ is not even a weak local minimum in problem \eqref{GrindEQ__5_5_}. 

Consequently, by virtue of Example 5.2 we come to the conclusion that the statements of Theorem 2.1 in the paper \cite{3} were strengthened in the forms \eqref{GrindEQ__3_16_}, \eqref{GrindEQ__3_22_} and \eqref{GrindEQ__4_14_}, while the statements of Theorem 2.3 in the paper \cite{3} were strengthened by means of Theorem 3.6 and part (jj) of Theorem  3.7.

It is important to note that in problem \eqref{GrindEQ__5_5_} along the extremal $\bar{x}\left(t\right)=0,\, \, t\in \left[0,\, 1\right]$, the expressions in the forms \eqref{GrindEQ__3_1_}, \eqref{GrindEQ__3_2_} and \eqref{GrindEQ__3_5_} are not fulfilled for any point $\left(\theta ,\, \bar{\lambda },\, \eta \right)$ of the set $\left[0,1\right]\times \left(0,\, 1\right)\times R^{2} \backslash \left\{0\right\}$, since $\frac{\bar{\lambda }}{\bar{\lambda }-1} \left(1-\frac{\bar{\lambda }}{\bar{\lambda }-1} \right)\ne 0$, for all $\bar{\lambda }\in \left(0,\, 1\right)$. Therefore, by virtue of problem 5.5 we confirm that degeneration in the form \eqref{GrindEQ__4_13_}, and also degenerations in the forms  \eqref{GrindEQ__3_17_}, \eqref{GrindEQ__3_18_} and \eqref{GrindEQ__3_21_}, generally speaking, do not yield the degenerations in the forms \eqref{GrindEQ__4_1_}, \eqref{GrindEQ__3_1_}, \eqref{GrindEQ__3_2_} and \eqref{GrindEQ__3_5_}, respectively. 

 \textbf{Example 5.3. }Let us consider the problem 
\begin{equation} \label{GrindEQ__5_7_} 
J\left(x\left(\cdot \right)\right)=\int _{0}^{1}\left(\left(1-t\right)\dot{x}^{3} -3x\right)dt\to {\mathop{\min }\limits_{x\left(\cdot \right)}} ,\, \, \, \, \, x\left(0\right)=0, \, \, \, x\left(1\right)=1,    
\end{equation} 
where $x\left(\cdot \right)\in KC^{1} \left(\left[0,\, 1\right],\, R\right),\, L\left(\cdot \right)=\left(1-t\right)\dot{x}^{3} -3x$.

 Obviously, the admissible function $\bar{x}\left(t\right)=t$, $t\in \left[0,\, 1\right]$ is an extremal in problem \eqref{GrindEQ__5_7_}. It is a weak local minimum. This follows from the fact that for an arbitrary increment $\Delta x\left(\cdot \right)\in KC^{1} \, \left(\, \left[0,\, 1\right],\, R\right)$ for which $\Delta x\left(0\right)=\Delta x\left(1\right)=0$ and $\left\| \Delta \, \dot{x}\left(t\right)\right\| _{L_{\infty } \left(\left[0,\, 1\right],\, R\right)} \le 3$, the inequality
\[\Delta J\left(\bar{x}\left(\cdot \right)\right)=\int _{0}^{1}\left(\Delta \dot{x}\left(t\right)\right)^{2} \left(3+\Delta \, \dot{x}\left(t\right)\right)dt\ge 0\] 
is fulfilled.

\noindent Let us verify the validity of Theorem 2.4 proved in \cite{3}. Since $E\left(\bar{L}\right)\left(t,\, \xi \right)=\left(1-t\right)\, \xi ^{2} \left(\xi +3\right)$ and $\bar{L}_{\dot{x}\dot{x}} \left(t\right)=6\left(1-t\right)$, then at $t=1-0$ there are equalities (2.5) and (2.6) in the paper \cite{3}. Therefore, it is clear that for $\bar{x}\left(t\right)=t$, $t\in \left[0,\, 1\right]$, at the point $t=1$ and for all  $\xi \in R$ all the assumptions of Theorem 2.4 in the paper \cite{3} are fulfilled. Then, by virtue of the statement of Theorem 2.4 in the paper \cite{3}  there exists such $\delta >0$ that the inequality$-\xi ^{2} \left(\xi +6\right)\ge 0$ is valid for all $\xi \in \left(-\delta ,\, \delta \right)$. Obviously, there is no such a number $\delta >0$ that the last inequality would be fulfilled. However, the statement of part (j) of Theorem 3.7, i.e. inequality \eqref{GrindEQ__3_20_} is fulfilled for example for $\delta =6$. Consequently, the statement formulated for $\tau _{-} $ in Theorem 2.4 of the paper \cite{3} is wrong; furthermore, Theorem 3.6 and part (j) of Theorem 3.7 are correct and strengthened version of Theorem 2.4 in the paper \cite{3}. In a similar way we confirm that parts (i) and (ii) of Theorem 3.3 are correct and strengthened version of Theorem 2.2 in the paper \cite{3}. 

 \textbf{Example 5.4. }Let us consider the problem
\begin{equation} \label{GrindEQ__5_8_} 
J\left(x\left(\cdot \right)\right)=\int _{0}^{1}\left(\left(\dot{x}_{1} -\dot{x}_{2}^{3} \right)^{2} +x_{1} \dot{x}_{2}^{2} \right)\, dt\to {\mathop{\min }\limits_{x\left(\cdot \right)}} ,\, \, \, \, x_{i} \left(0\right)=x_{i} \left(1\right)=0, \, \, \, i=1,\, 2,  
\end{equation} 
where $x=\left(x_{1} ,\, x_{2} \right)^{T} ,\, \, x\left(t\right)\in KC^{1} \left(\, \left[0,\, 1\right],\, R^{2} \right),\, \, L\left(\cdot \right)=\left(\dot{x}_{1} -\dot{x}_{2}^{3} \right)^{2} +x_{1} \dot{x}_{2}^{2} $.

 We study the admissible function $\bar{x}\left(t\right)=0,\, \, t\in \left[0,\, 1\right]$ for a minimum. Along this function allowing for \eqref{GrindEQ__1_6_} and \eqref{GrindEQ__2_10_} we have
\[\bar{L}\left(t\right)=0,\, \bar{L}_{x} \left(t\right)=\bar{L}_{\dot{x}} \left(t\right)=0,\, \, L_{x} \left(t,\, \bar{x}\left(t\right),\, \dot{\bar{x}}\left(t\right)+\xi \right)=\left(\xi _{2}^{2} ,\, 0\right)^{T} ,\] 
\begin{equation} \label{GrindEQ__5_9_} 
E\left(\bar{L}\right)\, \left(t,\, \xi \right)=\left(\xi _{1} -\xi _{2}^{3} \right)^{2} ,\, \, \, \xi ^{T} L_{\dot{x}\dot{x}} \left(t\right)\xi =2\xi _{1}^{2} ,\, \, \Delta \bar{L}_{x} \left(t,\, \xi \right)=\left(\xi _{2}^{2} ,\, 0\right)^{T} .   
\end{equation} 

Obviously, the admissible function $\bar{x}\left(\cdot \right)=0$ is an extremal in problem \eqref{GrindEQ__5_8_} and along it the Weierstrass condition \eqref{GrindEQ__1_5_} is fulfilled.

 Considering \eqref{GrindEQ__5_9_}, we easily get that along the extremal $\bar{x}\left(\cdot \right)=0$ in problem \eqref{GrindEQ__5_8_} for $\bar{\lambda }=\frac{1}{2} $ and for all points $\left(t,\, \eta \right)\in \left(0,\, 1\right)\times \, \left\{\, \, \left(\eta _{2}^{3} ,\, \eta _{2} \right)^{T} :\eta _{2} \in R\, \right\}\, $ the degeneration in the form \eqref{GrindEQ__4_1_} is fulfilled. However, the degeneration in the form \eqref{GrindEQ__4_13_} is fulfilled only for $\eta =0$. Consequently, unlike Example 5.2 we confirm that, generally speaking, the degenerations in the forms \eqref{GrindEQ__3_1_}, \eqref{GrindEQ__3_2_}, \eqref{GrindEQ__3_5_} and \eqref{GrindEQ__4_1_} do not imply the degenerations in the forms \eqref{GrindEQ__3_17_}, \eqref{GrindEQ__3_18_}, \eqref{GrindEQ__3_21_} and \eqref{GrindEQ__4_13_}, respectively. 

 We continue our study. Since the degeneration in the form \eqref{GrindEQ__4_13_} is fulfilled only for $\eta =0$, then Theorem 4.3 keeps the extremal $\bar{x}\left(\cdot \right)=0$ among the contenders for a minimum. Note that Theorem 4.2 is also ineffective while studying the extremal $\bar{x}\left(\cdot \right)=0$.

 We now apply Theorem 4.1. Since along the extremal $\bar{x}\left(\cdot \right)=0$ for$\bar{\lambda }=\frac{1}{2} $ and for all $\left(t,\, \eta \right)\in \left(0,\, 1\right)\times \left\{\, \left(\eta _{2}^{3} ,\, \eta _{2} \right)^{T} :\eta _{2} \in R\, \right\}$ the degeneration in the form \eqref{GrindEQ__4_1_} is fulfilled, then allowing for \eqref{GrindEQ__5_9_} the necessary condition \eqref{GrindEQ__4_2_} takes the form: $\eta _{2}^{5} =0$ for all $\eta _{2} \in R$. Consequently, necessary minimum condition \eqref{GrindEQ__4_2_} is not fulfilled for all $\eta_2 \ne 0$, so, by virtue of Theorem 4.1 the extremal $\bar{x}\left(\cdot \right)=0$ is not even a weak local minimum in problem \eqref{GrindEQ__5_8_}.

In conclusion, we consider it promising to obtain analogues of some results of this paper, for example, the analogues of Theorems 4.1-4.3 in the theory of singular optimal controls.

\end{document}